\documentclass[10pt]{article}
\usepackage[T1]{fontenc}
\usepackage[nohead,margin=1.0in]{geometry}
\usepackage{graphicx,caption,subcaption}

\usepackage{amsmath,amssymb,amsthm, amsfonts}
\usepackage{url}
\usepackage{multicol}
\usepackage{multirow}
\usepackage{parskip}
\usepackage[sorting=none,sortcites,maxbibnames=99]{biblatex}
\usepackage{xcolor}
\usepackage{color}
\usepackage{nicefrac}
\usepackage{hyperref}
\usepackage{mathrsfs}
\usepackage{siunitx}
\usepackage{makecell}

\usepackage{algorithmic}
\usepackage{algorithm}
\usepackage{array}

\usepackage{tikz}
\usetikzlibrary{tikzmark}
\usepackage{pgfplots}
\pgfplotsset{compat=1.18}

\newcommand{\N}{\mathbb{N}}

\newcommand{\R}{\mathbb{R}}
\newcommand{\C}{\mathbb{C}}
\renewcommand{\d}{\, \text{d}}

\theoremstyle{plain}

\numberwithin{equation}{section}

\DeclareSIUnit{\vacuumpermeability}{\text{$\mu_0$}}

\DeclareMathOperator*{\argmax}{arg\,max}
\DeclareMathOperator*{\argmin}{arg\,min}

\begin{document}

\title{Learned Discrepancy Reconstruction and Benchmark Dataset for Magnetic Particle Imaging}

\author{Meira Iske\thanks{Center for Industrial Mathematics, University of Bremen, Bremen, Germany.  Corresponding author: \texttt{iskem@uni-bremen.de}} \and Hannes Albers\footnotemark[1] \and Tobias Knopp\thanks{Section for Biomedical Imaging, University Medical Center Hamburg-Eppendorf, Hamburg, Germany.}\textsuperscript{ ,}\thanks{Institute for Biomedical Imaging, Hamburg University of Technology, Hamburg, Germany.}\textsuperscript{ ,}\thanks{Fraunhofer Research Institution for Individualized and Cell-Based Medical Engineering IMTE, Lübeck, Germany.} \and Tobias Kluth\footnotemark[1]}

\date{}

\maketitle

\begin{abstract}
Magnetic Particle Imaging (MPI) is an emerging imaging modality based on the magnetic response of superparamagnetic iron oxide nanoparticles to achieve high-resolution and real-time imaging without harmful radiation. One key challenge in the MPI image reconstruction task arises from its underlying noise model, which does not fulfill the implicit Gaussian assumptions that are made when applying traditional reconstruction approaches. To address this challenge, we introduce the Learned Discrepancy Approach, a novel learning-based reconstruction method for inverse problems that includes a learned discrepancy function. It enhances traditional techniques by incorporating an invertible neural network to explicitly model problem-specific noise distributions. This approach does not rely on implicit Gaussian noise assumptions, making it especially suited to handle the sophisticated noise model in MPI and also applicable to other inverse problems. 
To further advance MPI reconstruction techniques, we introduce the MPI-MNIST dataset — a large collection of simulated MPI measurements derived from the MNIST dataset of handwritten digits. The dataset includes noise-perturbed measurements generated from state-of-the-art model-based system matrices and measurements of a preclinical MPI scanner device. This provides a realistic and flexible environment for algorithm testing.  
Validated against the MPI-MNIST dataset, our method demonstrates significant improvements in reconstruction quality in terms of structural similarity when compared to classical reconstruction techniques.
\end{abstract}

\section{Introduction}
Magnetic Particle Imaging (MPI) is a promising imaging modality that exploits the magnetization behaviour of super paramagnetic nanoparticles~\cite{Gleich2005}. The nonlinear magnetization response of the particles to externally applied magnetic fields induces a voltage, which is measured by receive coils. The acquired signal can be translated to the spatial distribution of the particle concentration, thereby enabling image reconstruction~\cite{Knopp2017}.
MPI offers high spatial resolution without the use of harmful radiation. Additionally, it exhibits high sensitivity, with the signal remaining unaffected by background noise from surrounding tissues~\cite{Buzug2012}. Due to its rapid data acquisition capabilities, this imaging technique enables real-time imaging~\cite{Gleich2008,Knopp2016_2}. These characteristics make MPI advantageous for a range of clinical applications, such as vascular applications~\cite{Herz2018} including blood flow imaging~\cite{Weizenecker2009, kaul2018magnetic} and instrument tracking~\cite{Haegele2012}. For instance, it has applications in the precise positioning of stents~\cite{Herz2019} or catheters~\cite{Rahmer2017, ahlborg2022first}.  
Additionally, by labeling cells with super paramagnetic particles, MPI can be applied for cell tracking~\cite{Bulte2015,Sehl2020,bulte2022vivo}, which is particularly relevant for cancer detection and cancer therapy~\cite{Yu20217,bauer2016high}. 

The reconstruction of the particle concentration requires solving a linear ill-posed inverse problem~\cite{Knopp2010, Knopp2008} by the means of the system matrix, which constitutes the forward operator. However, appropriate reconstruction quality is concerned with major aspects of the modeling chain such as a careful determination of the system matrix and accurate modeling of noise conditions. 
The former aspect can be accomplished through a data-driven calibration approach. Although this method accounts for all relevant physical properties, it introduces noise and its measurement process is time-consuming. These limitations have led to the development of model-based system matrices, that determine the MPI system matrix based on physical models. The latter allows for noise-free system matrices and potentially much faster generation times.
Initial model-based methods~\cite{Knopp2010_model, Knopp2010_model2D} employed the simplified equilibrium model, which assumes no relaxation time for the magnetic nanoparticles. More advanced approaches include the effect of Brownian~\cite{Yoshida2012,Kluth2018} and Néel rotation~\cite{Kluth2019_towards, Albers2022} of the magnetic nanoparticles on its magnetic moment by solving the Fokker-Planck equation. For additional details on model-based system matrices, we refer to \cite{Kluth2018,maass2024} and references therein. 

The second essential aspect of the modeling chain focuses on the reduction of unknown noise introduced by the MPI scanner, which exhibits characteristics that cannot generally be described by a standard normal distribution. Previous approaches~\cite{Them2015, Knopp2019_bg} have addressed this component by applying a background correction from scanner measurements prior to the image reconstruction procedure. 
The conventional reconstruction objective involves minimizing the Tikhonov functional, typically achieved through the application of the regularized Kaczmarz algorithm~\cite{Weizenecker2009,Knopp2016_2,Kluth2019_towards}.
Subsequent research~\cite{Knopp2010, Kluth2019} 
adapts the formulation of the common objective by introducing a linear transform to the data discrepancy term of the Tikhonov-type framework. 
While the linear transform in~\cite{Knopp2010} is defined by the row-energy of the system matrix, \cite{Kluth2019} expands upon this by applying a whitening operation that directly includes the scanner noise characteristics in a data-based fashion.
A derivation of the Tikhonov objective from a Bayesian perspective by estimating the maximum a posteriori~(MAP) amounts to making a standard Gaussian assumption on the noise distribution. Similarly, for the adapted objective in~\cite{Knopp2010, Kluth2019} one assumes Gaussian distributed noise with a covariance matrix depending on the transformation matrix. 

Beyond noise-specific methods, further reconstruction approaches address system matrix deviations and data structure priors. For instance, total least squares methods~\cite{Kluth2017} and the RESESOP algorithm~\cite{Blanke2020, Nitzsche2022} focus on system matrix inaccuracies. Prior data knowledge can be introduced by applying a TV penalty term~\cite{Bathke2017, Storath2017} or by using an $\ell^1$-data discrepancy term~\cite{Kluth2020} in the formulation of the objective. In recent years, there has been significant interest in the development of learning-based methods for MPI, particularly those leveraging deep learning techniques. These data-driven methods seem extremely promising due to their ability of capturing problem specific data structures and have been applied by following various strategies. Post-processing tools utilize neural networks to refine images reconstructed by traditional methods~\cite{Knopp2023}. 
End-to-end learning approaches have gained attention for their ability to directly learn the inverse of the system matrix~\cite{Gladiss2022,Gladiss2022_1,Koch2019}, where plug-and-play methods integrate neural networks into an iterative reconstruction framework~\cite{Askin2022,Gungor2023,Tsanda2024}. 
In this context, the neural network is typically responsible for denoising in image space. Similarly, deep equilibrium approaches apply iterative schemes as in~\cite{gungor2023deq,Gungor2024}, by employing a fixed-point model where the output of the neural network represents the equilibrium point of the iterative process. The unsupervised Deep Image Prior (DIP) approach~\cite{Dittmer2021} has also been applied in MPI, where the neural network itself acts as a prior, capturing image structures without the need for extensive training data.

The application of suitable reconstruction schemes is crucial in order to obtain reasonable reconstructions. The development of effective reconstruction techniques, in turn, relies on the availability of validation data. However, the limited number of MPI scanners poses a challenge in providing the research community with sufficient access to reference data. So far, the publicly available OpenMPI~\cite{Knopp2020} dataset has been widely used for the evaluation of reconstruction algorithms~\cite{Kluth2019, Dittmer2021, Kluth2020, Shang2022}. 
Despite representing an initial step toward establishing a data foundation for the evaluation of image reconstruction approaches, the number of available MPI data samples 
is still limited. Moreover, there is a general lack of paired ground truth measurements in MPI. This presents a significant challenge when it comes to the conception of learning-based approaches, typically requiring a large amount of data.

In this work we revisit two key challenges in MPI reconstruction.  First, we propose a data-based reconstruction approach that directly accounts for complex and problem-specific noise models, thereby addressing one aspect of the modeling chain. Unlike previous methods that use simplified assumptions on the noise model, leading to standard discrepancy terms within the Tikhonov framework, our approach incorporates a learned discrepancy term. This term relies on an invertible neural network trained to learn the unknown noise distribution. This results in an objective consistent with the variational regularization framework, which is tailored to the underlying noise characteristics that were captured during training. While related learning-based approaches in other contexts, such as~\cite{Pinetz2021, Calatroni2017}, have focused on learning a parametrization of discrepancy functions for selected noise distribution classes, our method derives the discrepancy functional from a Bayesian perspective. This integration of the learned noise distribution into the likelihood term enables a reconstruction method that extends beyond Gaussian noise assumptions and is transferable to a wide range of inverse problems.
Second, we introduce the MPI-MNIST dataset, based on the MNIST dataset of handwritten digits~\cite{lecun1998mnist}. 
It includes both simulated MPI measurements obtained by using a model-based system matrix in combination with real scanner noise from the preclinical MPI system~(Bruker, Ettlingen, Germany). This enables a realistic yet verifiable testing environment for reconstruction algorithms, where the contained noise measurements serve as a basis for our proposed data-based reconstruction method. With the MPI-MNIST dataset we intend to make a further step toward closing the existing data gap in MPI.

This work is structured as follows: 
In Section~\ref{sec:methods} we cover different reconstruction methods, where Section~\ref{subsec:standard_rec} focuses on common approaches in the context of MPI image reconstruction. In Section~\ref{sec:lmap} we introduce the Learned Discrepancy Approach. Section~\ref{sec:physical_model} details the generation of the model-based system matrix, which is followed by the process of measuring noise samples in Section~\ref{subsec:NoiseMeas} that both build the foundation of the MPI-MNIST dataset. The dataset’s structure is thoroughly described in Section~\ref{sec:dataset}, which is applied to our reconstruction approach in Section~\ref{sec:num}. In this section, we present a qualitative and quantitative comparison of our approach with the methods from Section~\ref{subsec:standard_rec}.

\section{Methods}\label{sec:methods}

We consider reconstruction methods with broad applicability to linear inverse problems, where ${X \subseteq R^N}$, $ {Y \subseteq \R^M}$ and ${A: X \rightarrow Y}$ a linear operator satisfying
\begin{equation}\label{eq:inv_prob}
    Ax = y
\end{equation}
for $x \in X$ and $y \in Y$. In this setting, we aim to reconstruct $x$ from noise-perturbed data $y^\delta = y + \eta$ with $\| \eta \| \leq \delta$. 

\subsection{Standard Reconstruction Approaches}\label{subsec:standard_rec}

Typically, inverse problems of type~\eqref{eq:inv_prob} are solved by the concept of variational regularization, namely by minimizing 
\begin{equation*}
    \mathcal{T}(x) = \mathcal{D}(Ax,y^\delta) + \mathcal{R}(x),
\end{equation*}
where $\mathcal{D}: Y \times Y \longrightarrow [0,\infty]$ measures the data discrepancy and $\mathcal{R}: X \longrightarrow (-\infty,\infty] $ denotes the regularization function. Specific choices for both functions lead to the definition of the \emph{Tikhonov functional}
\begin{equation}\label{eq:tikh_obj}
    \mathcal{T}_\alpha^{y^\delta}(x) = \| y^\delta - Ax\|^2_Y + \alpha \| x\|^2_X,
\end{equation}
with $\alpha > 0$. 
The analytical solution of~\eqref{eq:tikh_obj} is obtained by 
\begin{equation*}
    \hat{x}_\mathrm{Tikhonov} = (A^\mathsf{H} A + \alpha I)^{-1}A^\mathsf{H} y^{\delta},
\end{equation*}
where $A^\mathsf{H}$ denotes the conjugate transpose of $A$. 

A frequently applied reconstruction approach in MPI image reconstruction is the so-called 
\emph{regularized Kaczmarz}~(RK) algorithm, as proposed in~\cite{Dax1993}. This well-established reconstruction method iteratively converges to the minimizer of~\eqref{eq:tikh_obj} and
is often chosen due to its rapid convergence while still delivering reconstructions of notably high quality as demonstrated in~\cite{Kluth2019_towards, Kluth2019}. 

Beyond the conventional regularized Kaczmarz algorithm, \cite{Kluth2019} have shown the success of the whitened variant of this algorithm. Its objective is expressed as
\begin{equation}\label{eq:wrk_obj}
    \hat{x} = \argmin_x \frac{1}{2} \| W(y^\delta - Ax) \|^2 + \frac{\alpha}{2} \| x\|^2,
\end{equation}
where $W \in \R^{K \times K}$ a diagonal matrix with entries $w_{j,j} = \min_k \nicefrac{\mathrm{std}(\eta_k)}{\mathrm{std}(\eta_j)}$ for $j = 1,\hdots,K$. Here, $\mathrm{std}(n_j)$ is computed for each noise component $j$ over the additive measurement noise on the train set. 
The minimizer in~\eqref{eq:wrk_obj} is again determined via the regularized Kaczmarz algorithm, by defining $A_W := WA$ and $y_W^\delta = Wy^\delta$ in~\eqref{eq:tikh_obj}. 
This approach balances the impact of different noise features within the optimization process. Following the terminology of~\cite{Kluth2019}, we refer to the latter reconstruction method as \emph{whitened regularized Kaczmarz}~(WRK).

\subsection{Learned Discrepancy Approach}\label{sec:lmap}

\noindent In this section, we introduce the Learned Discrepancy Approach~(LDA), a learning-based reconstruction method, which is broadly applicable to inverse problems following an operator equation as in~\eqref{eq:inv_prob}.
The proposed method is derived from a Bayesian perspective and 
aims to encode the intrinsic noise distribution of the underlying measurement noise. Since MPI noise characteristics can at this point only be inferred through measurement data, the approach seems especially promising for MPI image reconstruction tasks. It is particularly effective with the data collection contained in the MPI-MNIST dataset as proposed in Section~\ref{sec:dataset}.
By leveraging MPI-specific noise characteristics, our approach holds significant promise for enhancing reconstruction accuracy. 

We assume that $\eta$ follows a probability density function~(pdf) $p_H: Y \longrightarrow \R_0^+$ and $x$ follows the prior pdf $p_X: X \longrightarrow \R_0^+$. In the following, we will further assume that the random variables $x \sim p_X$ and $\eta \sim p_H$ are stochastically independent. 

When performing reconstruction tasks, from a probabilistic point of view, we are interested in the posterior distribution $p(x|y^\delta)$ in order to identify the value of $x$ that maximizes the probability given the observed data $y^\delta$. According to Bayes' Theorem, the posterior distribution is given by
\begin{equation}\label{eq:bayes_thm}
    p(x|y^\delta) = \frac{p(y^\delta | x)p_X(x)}{p_Y(y^\delta)} \propto p(y^\delta | x) p_X(x), 
\end{equation}
where $p_Y$ describes the pdf of $y \in Y$. The maximum a posteriori~(MAP) estimate $\hat{x} \in X$ is obtained by maximizing the posterior distribution, i.e.,
\begin{align}\label{eq:max_posterior}
    \hat{x} &= \argmax_x \ p(x | y^\delta) = \argmax_x \ \log p(x | y^\delta) = \argmin_x \ - \log p(x | y^\delta).
\end{align}
Combining~\eqref{eq:bayes_thm} and~\eqref{eq:max_posterior} implies
\begin{equation}\label{eq:MAP}
    \hat{x} = \argmin_x \ - \log p(y^\delta | x) - \log p_X(x).
\end{equation}
Due to the independence of $x$ and $\eta$, the likelihood can be expressed as $p(y^\delta | x) = p_H(y^\delta -Ax)$. Commonly, $\eta$ is assumed to follow a Gaussian distribution $p_H = \mathcal{N}(0,I)$, which, together with~\eqref{eq:MAP} leads to the objective
\begin{equation}\label{eq:tikh_obj2}
    \hat{x} = \argmin_x \  \frac{1}{2}\| y^\delta - Ax\|^2 + \mathcal{R}(x),
\end{equation}
as known from the Tikhonov regularization with regularizing functional $\mathcal{R}(x) =- \log p_X(x)$. By making a Gaussian assumption on the prior, i.e., assuming $p_X = \mathcal{N}(0,\frac{1}{\alpha}I)$ for $\alpha > 0$, we obtain $\mathcal{R}(x) = \frac{\alpha}{2} \| x\|^2$ as in~\eqref{eq:tikh_obj}.

However, in practical scenarios, where the noise follows some sophisticated structures due to given physical circumstances or external influences, it can be restrictive to assume Gaussian distributed noise. Our aim is to reformulate and generalize~\eqref{eq:tikh_obj} by allowing for arbitrary pdfs $p_H$. By exploiting the features of invertible neural networks~(INNs), which can learn approximations of point-wise evaluations of an unknown $p_H$, we adapt the discrepancy term in~\eqref{eq:tikh_obj}, which is tailored to the noise distribution $p_H$. 

 Let $z \in \R^M$ be a random variable following a known and tractable probability density function $p_Z$, which we will refer to as the latent distribution. Let $g: \R^M \longrightarrow \R^M$ be a diffeomorphism, i.e., an invertible and differentiable mapping. Furthermore, let $g$ be such that $z = g(\eta)$. By the change-of-variable theorem it holds 
\begin{equation}\label{eq:cov}
    p_H(\eta) = p_Z(g(\eta)) | \det J_{g}(\eta)|,
\end{equation}
where $J_{g}(\eta) := \frac{\partial g}{\partial \eta} $ denotes the Jacobian of $g$. Therefore, the underlying probability density function $p_H$ can be expressed as a transformation of the latent distribution $p_Z$ through the diffeomorphism $g$.

To approximate the true underlying distribution $p_H$, we train an invertible neural network $\varphi_\theta: \R^M \longrightarrow \R^M$ parameterized by $\theta \in \Theta$, where $\Theta$ denotes the parameter space of the neural network in order to model the diffeomorphism $g$ given in~\eqref{eq:cov}. The goal is to find a transformation $\varphi_\theta$ that maps an input $\eta$ to a latent variable $z$ such that $p_\theta(\eta) \approx p_H(\eta)$ can be approximated, where $p_\theta(\eta)$ is the model distribution defined by
\begin{equation*}
    p_\theta(\eta) = p_Z(\varphi_\theta(\eta)) |\det(J_{\varphi_\theta}(\eta)) |.
\end{equation*}

The network is trained on a dataset $\{ \eta^{(i)}\}_{i=1}^L$ by maximizing the log-likelihood
\begin{align}\label{eq:loss}
    &\max_{\theta \in \Theta} \mathcal{L}(\theta) \\  \text{ with } \mathcal{L}(\theta) &=  \sum_{i=1}^L \log(p_\theta(\eta^{(i)})) = \sum_{i=1}^L \left( \log p_Z(\varphi_\theta(\eta^{(i)})) + \log | \det J_{\varphi_\theta}(\eta^{(i)}) |\right), \notag
\end{align}
which is equivalent to minimizing the Kullback-Leibler divergence $D_\mathrm{KL}(p_H \| p_\theta)$.

By choosing our latent distribution $p_Z$ to be of Gaussian nature, together with~\eqref{eq:MAP} we have 
\begin{align*}
    \hat{x} &= \argmin_x \ -\log p_H(y^\delta-Ax)-\log p_X(x) \\
            &= \argmin_x \ -\log p_Z(\varphi_\theta(y^\delta-Ax)) -\log | \det J_{\varphi_\theta}(y^\delta-Ax)| - \log p_X(x) \\
            &= \argmin_x \ \frac{1}{2} \| \varphi_\theta(y^\delta -Ax) \|^2 -\log | \det J_{\varphi_\theta}(y^\delta-Ax)|  - \log p_X(x),
\end{align*}
which is a modified version of the typical MAP estimation that incorporates a normalizing flow to account for a more complex noise distribution beyond the standard Gaussian assumption. We define 
\begin{equation*}
   \mathcal{D}_\theta(Ax, y^\delta) := \frac{1}{2} \| \varphi_\theta(y^\delta -Ax) \|^2 -\log | \det J_{\varphi_\theta}(y^\delta-Ax)|, 
\end{equation*} so that our objectives writes as
\begin{equation}\label{eq:lmap_obj}
    x_\mathrm{LDA} = \argmin_x \ \mathcal{D}_\theta(Ax, y^\delta) - \mathcal{R}(x),
\end{equation}
where $x_\mathrm{LDA}$ is the minimizer of the Learned Discrepancy Approach~(LDA). The LDA incorporates a learned discrepancy functional $\mathcal{D}_\theta$ and extends the objective in~\eqref{eq:tikh_obj} by the log-determinant of the trained network and an evaluation of the network within the discrepancy term. Therefore, we refer to $\varphi_\theta$ as the \emph{discrepancy network}. 

\section{Physical Model}\label{sec:physical_model}

\noindent As a preparatory step for data generation and for testing the proposed algorithm in Section~\ref{sec:lmap} in the context of MPI, 
we simulate system matrices in a number of different settings. This simulation concept is in line with the one given in~\cite{Kluth2019_towards}, which we extend by using a linear combination of simulations in order to obtain a better approximation of a measured system matrix. One such system matrix was simulated to provide, to the best of our knowledge, the closest approximation to a measured system matrix in the same setting. We therefore treat this reference simulation as the best possible approximation to a noise-free system matrix. The remaining simulated system matrices are used for numerical experiments, where we assume that only a system matrix of lower quality than the real-life physical system is available. This allows for testing in contexts such as operator correction or post-processing of artifacts resulting from incomplete or inaccurate modeling.

\subsection{System Matrix Simulation Framework} \label{subsec:SM}

\noindent The voltage induced in a receive coil of an MPI scanner with sensitivity profile $p: \R^3 \longrightarrow \R^3$ can be described by
\begin{equation*}
    \Tilde{v}(t) = -\mu_0 \int_\Omega c(x) p(x)^\mathsf{T} \frac{\partial}{\partial t}\bar{m}(x,t) \ \d x
\end{equation*}
in $\mathrm{V}$, where $c: \Omega \rightarrow \R_0^+ $ in $\mathrm{mol/L}$ is the concentration of the magnetic nanoparticles, $\Omega \subset \R^3$ the imaging volume and $\bar{m}: \R^3 \times [0,T]\rightarrow \R^3$ the mean magnetic moment in $10^{-3}\mathrm{Am^2 mol^{-1}}$. The direct feedthrough of the signal is filtered by an analog filter that can be expressed by a filter kernel $a:\R \rightarrow \R$, where $v = a \ast \Tilde{v}$ the filtered signal. The $T$-periodic signal $v$ can consequently be expanded into a Fourier series with coefficients

\begin{equation*}
    \hat{v}_k = -a_k\frac{\mu_0}{T}\int_\Omega c(x) p(x)^\mathsf{T} \int_0^T \left( \frac{\partial }{\partial t} \bar{m}(x,t)\right)\exp\left( \frac{-2\pi i t k}{T}\right) \ \d t \ \d x 
\end{equation*}

for $k \in \N_0$. By defining the system function 

\begin{equation*}
    s_k(x) = -a_k\frac{\mu_0}{T} \int_0^T p(x)^\mathsf{T}  \left( \frac{\partial }{\partial t} \bar{m}(x,t)\right)\exp\left( \frac{-2\pi i t k}{T}\right) \d t,
\end{equation*} 

the Fourier coefficients write as $\hat{v}_k = \int_\Omega c(x)s_k(x)\ dx$. Sampling from the system function then yields the system matrix $A = (s_k(x_\ell))_{k=0,\hdots,K-1;\ell=1,\hdots,N} \in \C^{K \times N}$. A crucial part of the computation of the components $s_k$ requires access to the mean magnetic moment, depending on the $T$-periodic external magnetic field $H_\mathrm{app}: \R^3 \times [0,T] \rightarrow \R^3$ in $\mathrm{Am^{-1}}$. The applied magnetic field combines the selection field $H_S: \R^3 \rightarrow \R^3$ with the drive field $H_D: \R^3 \times [0,T] \rightarrow \R^3$ by $H_\mathrm{app}(x,t) = H_S(x) + H_D(x,t)$. The magnetic mean writes as
\begin{equation}
    \bar{m}(x,t) = m_0 \int_{S^2}m f(m,x,t) \ \d m, 
\end{equation}
where $m_0 = V_C M_S$ with $V_C = \frac{1}{6}\pi D_\mathrm{core}^2$ and $f$ the solution to the Fokker-Planck equation 
\begin{equation*}
    \frac{\partial}{\partial t }f = \mathrm{div}_{S^2}\left( \frac{1}{2\tau} \Delta_{S^2}f\right) - \mathrm{div}_{S^2}(bf),
\end{equation*}
with $\tau > 0 $ the relaxation time constant and (velocity) field $b: S^2 \times \R^3 \times \R^2 \rightarrow \R^3$, where
\begin{align*}
    b(m,H_\mathrm{app},n) &=  p_1 \times m + p_2(m \times H_\mathrm{app}) \times m + p_3 (n \cdot m ) n \times m + p_4(n \cdot m)(m\times n) \times m,
\end{align*}
 with $n\in S^2$ the easy axis of the magnetic particles. In this setting, we consider a pure Néel rotation of the magnetic particles including anisotropy, defining the physical constants $p_i$, $i=1,\hdots,4$ by $p_1 = \Tilde{\gamma} \mu_0$, $p_2 = \Tilde{\gamma}\alpha \mu_0$, $p_3 = p_4 = 2\Tilde{\gamma}\frac{K_\mathrm{anis}}{M_S}$ and $\tau = \frac{V_C M_S}{2 k_B T_B \Tilde{\gamma}\alpha}$ with $\Tilde{\gamma} = \frac{\gamma}{1+\alpha^2}$. To model fluid particles, the easy axis and the anisotropy constant are selected in analogy to Model B3 of~\cite{Kluth2019_towards}: The easy axis $n(x) = \frac{H_S(x)}{| H_S(x) |}$ is chosen in the direction of the selection field and the anisotropy constant $K_\mathrm{anis}(x) = K_\mathrm{anis}^\mathrm{max} \frac{| H_S(x) |}{h}$, with constant $K_\mathrm{anis}^\mathrm{max}$ and $h:=\max_{x\in \Omega} |H_S(x)|$, is chosen proportional to the selection field.

The system matrix simulation was performed by solving the Fokker-Planck equation for each pixel as described in \cite{Albers2022simulating}. Different system matrices were created for a number of different assumed maximum N\'{e}el anisotropy constants $K_\mathrm{anis}^\mathrm{max}$ and particle core diameters $D_\mathrm{core}$. In particular, system matrices were simulated for the following equidistant parameters:
\begin{align*}
    D_\mathrm{core} &\in \mathcal{D} = \{ 16,18,\dots, 24\,\si{nm}\}.\\
    K_\mathrm{anis}^\mathrm{max} &\in \mathcal{K} = \{ 500, 800, \dots, 3800\, \si{J/m^3}\}.
\end{align*}
We denote a system matrix corresponding to the combination of these parameters by $A_{\mathcal{D}(i),\mathcal{K}(j)}$ with $\mathcal{D}(i)$ the $i$-th element of $\mathcal{D}$ and $\mathcal{K}(j)$ the $j$-th element of $\mathcal{K}$.

To model the most accurate polydisperse system matrix, 
the system matrices corresponding to the different parameter combinations were equipped with weights $w_{i,j} \geq 0$, which were estimated in a least-squares sense by
\begin{align}\label{eq:lstsq}
    W^\ast = \argmin_{W} \| A_\mathrm{calib} - \sum_{i,j} w_{i,j} A_{\mathcal{D}(i),\mathcal{K}(j)} \|^2,
\end{align}
where $W = (w_{i,j})_{i=1,\hdots,5,j=1,\hdots,12}$ and $A_\mathrm{calib}$ is a suitable measured reference matrix. The resulting matrix $ A_\mathrm{opt} = \sum_{i,j} w^\ast_{i,j}A_{\mathcal{D}(i),\mathcal{K}(j)}$ with $W^\ast = (w^\ast_{i,j})_{i=1,\hdots,5,j=1,\hdots,12}$
will be treated as the most accurate available simulated matrix. In addition, the two system matrices with the lowest and with the highest anisotropy and diameter were included in the dataset to include mismatched or unrealistic system matrices. These are denoted $A_\mathrm{SP} := A_{\mathcal{D}(1),\mathcal{K}(1)}$ for the small parameters and $A_\mathrm{LP} := A_{\mathcal{D}(5),\mathcal{K}(12)}$ for the large parameters. Lastly, a standard equilibrium model was evaluated for $D_\mathrm{core} = 20 \, \si{nm}$ to serve as a crude approximation of the true system matrix, denoted by $A_\mathrm{eq}$.

Finally, all system matrices were scaled and time-shifted to fit $A_\mathrm{calib}$. The measured transfer function \cite{Thieben2024transfer} was added to each system matrix for its comparability to measurements. In total, we result in a set of system matrices where $A_\mathrm{opt}$ is as close as possible to the real-world situation, while the others contain simplifications or model errors of varying magnitude.

\subsection{Experimental Setup for $A_\mathrm{calib}$}

\noindent The preclinical MPI device (Bruker, Ettlingen, Germany), which was used for measuring $A_\mathrm{calib}$ has three send and three receive channels with excitation fields of $12\, \mathrm{mT}\mu_0^{-1}$ in both the $x$- and $y$-direction. The selected frequencies of the drive field are $f_x= 2.5/102 \,\mathrm{MHz}$ in the $x$-direction and  $f_y = 2.5/96 \,\mathrm{MHz}$ in $y$-direction. This results in a trajectory cycle of $652.8 \mu s$. The selection field gradient has a strength of $-1.0,-1.0$ and $2.0 \,\mathrm{T}\mu_0^{-1}\mathrm{m^{-1}}$ in $x$-, $y$-, and $z$-direction, respectively. The Lissajous trajectory is of cosine type, specifically cosine in $x$-direction and negative cosine in the $y$-direction. Signal reception uses the built-in 3D receive coils. The analog signal undergoes band-pass filtering and is digitized at a sampling rate of $2.5 \,\mathrm{MHz}$, where 1632 sampling points were captured for one Lissajous period, thus resulting in $K=817$ frequency components. A block average over 10 drive-field periods is applied during data acquisition. Spatial sampling of the system matrix takes place on a grid of size $17 \times 15$ with $N_\mathrm{coarse} = 255$ equally spaced measurement points, covering a $34 \times 34 \,\mathrm{mm}^2$ field of view. The used delta probe of size $2 \times 2 \times 2 \,\mathrm{mm}^3$ consists of Perimag. This tracer material contains an iron concentration of  $10 \,\mathrm{mg}_\mathrm{Fe}/\mathrm{mL}$. Altogether, this configuration yields the system matrix $A_\mathrm{calib} \in \C^{3 \times K \times N}$.

The simulation adapts to the experimental setup of the measured system matrix, while refining the spatial grid to a size of $85 \times 75$. The simulated system matrix was additionally subsampled for two grids of size $51 \times 45$ and $17 \times 15$, resulting in $A_\mathrm{opt}, A_\mathrm{SP}, A_\mathrm{LP}, A_\mathrm{eq} \in \C^{3 \times K \times N}$ with $N \in \{N_\mathrm{coarse} = 255, N_\mathrm{int} = 2295, N_\mathrm{fine} = 6375\}$.

\section{Noise Measurements} \label{subsec:NoiseMeas}
 \noindent To obtain a sufficient amount of noise measurements, the MPI scanner was left to run empty for a duration of approximately one hour while capturing the received voltage signals within this time period. The acquisition setup coincides with the setup of the reference system matrix as described in the previous section. In contrast to the measured system matrix, no block averaging was applied so that each individual acquisition within the Lissajous cycle resulted in a distinct noise frame $\bar{\eta}^{(j)} \in \C^{3\times K}$ for $j =1,\hdots,\num{4694400}$. Subsequent to the acquisition process, sequences of ten consecutive measurements were averaged by computing
\begin{equation*}
    \eta^{(i)} = \frac{1}{10}\sum_{j=1}^{10} \bar{\eta}^{(10(i-1)+j)} \in \C^{3 \times K} \quad \text{for } i = 1,\hdots,\num{469440}
\end{equation*}
to produce a single noise sample $\eta^{(i)}$. This is in line with the common measurement practice in MPI of performing multiple MPI measurements with subsequent averaging. The resulting noise samples are partitioned into five groups for the integration of the noise sampes to multiple parts of the dataset. Primarily, the captured samples serve as additive measurement noise of the simulated MPI measurements, as described in Section~\ref{subsubsec:sim_meas}. In addition, these samples can be used to simulate pixelwise deviations of the system matrix as in~\eqref{eq:H_SM}. Both additive measurement noise and system matrix deviations might require background corrections within the reconstruction process. Therefore, we further provide noise samples for background correction of both cases. To account for environmental factors such as temperature fluctuations during the measurement process, we employ alternating sequences of noise acquisitions for the system matrix and its background correction. Specifically, each sequence of noise measurements corresponding to a pixel in the system matrix is followed by an acquisition sequence of equal length reserved for background correction of the respective pixel including a short pause between the two sequences. Moreover, we provide a large amount of additional noise samples. These can be used e.g. for gaining statistical insights on the noise behaviour and are especially suitable for learning-based approaches as in Section~\ref{sec:lmap}. The organization of these five components is detailed in Section~\ref{sec:dataset}.

\section{MPI-MNIST: Dataset Generation and Structure}\label{sec:dataset}

\begin{figure}[t]
    \centering
    \includegraphics[width=\linewidth]{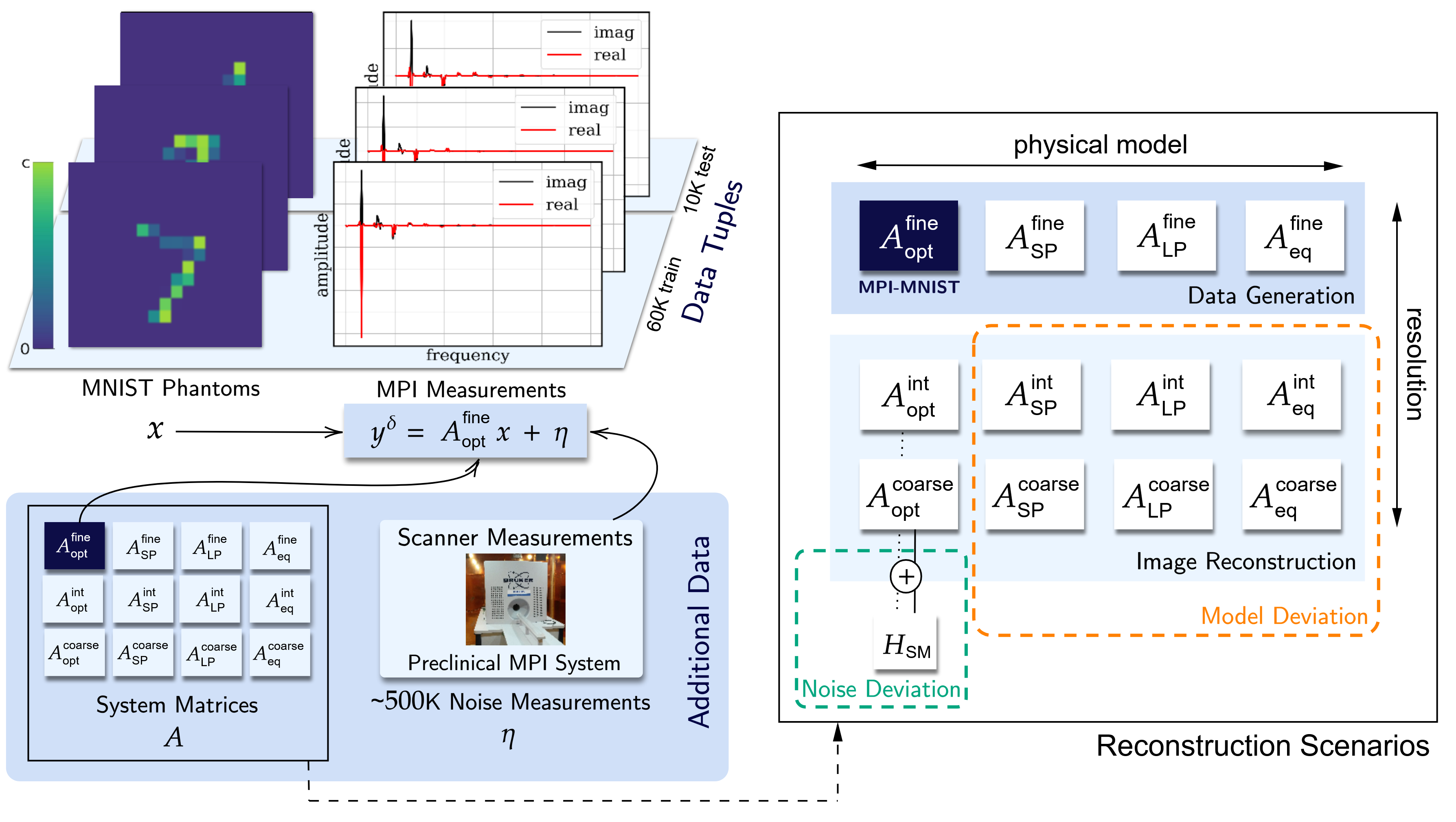}
    \caption{Overview of the MPI-MNIST dataset and included reconstruction settings.}
    \label{fig:sms}
\end{figure}

\noindent To address the existing data gap in MPI, we introduce the MPI-MNIST dataset containing data tuples of ground truth phantoms along with simulations of corresponding MPI measurements. An overview of the dataset is illustrated in Figure~\ref{fig:sms}. The measurement simulations result from the components described in Section~\ref{sec:physical_model} and Section~\ref{subsec:NoiseMeas}. The MPI-MNIST dataset comprises processed versions of all images included in the MNIST dataset of handwritten digits~\cite{lecun1998mnist}. 
The MNIST dataset itself contains \num{70000} labeled images of handwritten digits between zero and nine, thereby defining the size of the MPI-MNIST dataset. The images were chosen as a suitable ground truth dataset for MPI-MNIST due to the structural similarities of the images compared to common MPI phantoms such as blood vessels, particularly in terms of their sparsity. Consequently, we will refer to the two-dimensional MNIST images as \emph{phantoms} in the following discussion. The system matrices described in Section~\ref{sec:physical_model} are utilized for the simulation of MPI measurements in combination with the scanner data from Section~\ref{subsec:NoiseMeas}, which serve as measurement noise. Further system matrices underlying different physical models at various resolutions are provided. This increases the compatibility of the data with a broad range of reconstruction settings containing potential operator deviations, as illustrated on the right-hand side in Figure~\ref{fig:sms}. Additionally, different system matrices can be selected for data generation, as detailed in Appendix~\ref{subsec:custom_mpi-mnist}. The dataset further includes a large collection of noise data, separated into distinct components. These can serve for multiple applications, which are listed in Section~\ref{subsubsec:noise_meas}. 

The full dataset is publicly accessible via Zenodo~(\url{https://doi.org/10.5281/zenodo.12799417}). This iteration of the dataset features a specific setting, for which several parameters, e.g. particle concentration and resolution, had to be fixed. In 
this section we provide a comprehensive overview of this dataset's structure and data generation together with a comprehensive list of the provided files on the Zenodo website. 
We discuss potential modifications and extensions of the dataset in Appendix~\ref{subsec:custom_mpi-mnist}.

\subsection{Preprocessing of MNIST Phantoms}\label{subsubsec:prep_MNIST} 
The creation of data tuples in the MPI-MNIST dataset involves several adaptions to the original MNIST images, serving as the foundation for generating both noise-free and noise-perturbed simulated MPI measurements.
The MNIST phantoms are originally given in~\cite{lecun1998mnist} as images with a size of $28 \times 28 \, \mathrm{px}$. According to the adaption of the simulated system matrix to the measured reference system matrix in Section~\ref{sec:physical_model}, we aim to reconstruct on the same scale, i.e., reconstruct images of $17 \times 15 \, \mathrm{px}$ including an additional overscan of the system matrix. For this reason, the MNIST phantoms are first downsampled to the size of $11 \times 11 \, \mathrm{px}$ by applying nearest neighbour interpolation and subsequently framed by zero valued pixels. The resulting image is denoted by $x_\text{coarse}^{(i)} \in \R^{W_\mathrm{coarse} \times H_\mathrm{coarse}}$ with $W_\mathrm{coarse} = 17$, $H_\mathrm{coarse} = 15$ for $i=1,\hdots,\num{70000}$. The pixel values of each processed phantom are scaled to the desired pixel range $[0,c]$,  determined by the concentration $c = 10$ in $\text{mg}_\text{Fe}\text{mL}^{-1}$.
Therefore, $\big(x_\text{coarse}^{(i)}\big)_{j,k} \in [0,c]$ for $j=1,\hdots, W_\mathrm{coarse}, \ k = 1,\hdots,H_\mathrm{coarse}$. The samples $x_\text{coarse}^{(i)}$ will serve as ground truth images of our dataset. However, an upsampling of the phantoms is required for generation of the MPI data in order to circumvent an inverse crime~\cite{Kaipio2004}. We apply  nearest neighbour interpolation to $x_\text{coarse}^{(i)}$ by scaling each spatial dimension by a factor of $s = 5$. More precisely, we obtain $x_\mathrm{fine}^{(i)} \in \R^{W_\mathrm{fine} \times H_\mathrm{fine}}$ with $W_\mathrm{fine} = s \cdot W_\mathrm{coarse}, \ H_\mathrm{fine} = s \cdot H_\mathrm{coarse}$, where 
\begin{equation*}
    (x_\mathrm{fine}^{(i)})_{sj+u,sk+v} = (x_\mathrm{coarse}^{(i)})_{j,k} \qquad \text{for} \ u,v = 1,\hdots,s
\end{equation*}
and $j=1,\hdots,W_\mathrm{coarse}, \ k = 1,\hdots,H_\mathrm{coarse}$ defines the upsampling operation. In preparation for a matrix multiplication of the upsampled phantoms with the system matrix, we flatten the images to vector shape, resulting in $\bar{x}_\mathrm{fine}^{(i)} \in \R^{N_\text{fine}}$, $N_\mathrm{fine} = 6375$, where
\begin{equation*}
    (\bar{x}_\mathrm{fine}^{(i)})_\ell = (x_\mathrm{fine}^{(i)})_{j,k} \quad \text{with} \ \ell= j \cdot W_\mathrm{fine}+k, 
\end{equation*}
$j=1,\hdots, W_\mathrm{fine}$, $k=1,\hdots, H_\mathrm{fine}$. An analogous transformation is applied in order to obtain $\bar{x}^{(i)}_\mathrm{coarse}\in \R^{N_\mathrm{coarse}}$ from $x^{(i)}_\mathrm{coarse}$ with $N_\mathrm{coarse} = 255$.

\subsection{Simulation of MPI-MNIST Measurements}\label{subsubsec:sim_meas}

The data samples can be organized as tuples $(x^{(i)}, y^{(i)}, y^{\delta,(i)})$ for $i = 1,\hdots,\num{70000}$, where

\begin{itemize}
    \item[I] $x^{(i)} = \bar{x}_\text{coarse}^{(i)} \in \R^{N_\mathrm{coarse}}$ as in Section~\ref{subsubsec:prep_MNIST},
    \item[II] $y^{(i)} = [ A_1 \bar{x}^{(i)}_\mathrm{fine} \ | \ A_2 \bar{x}^{(i)}_\mathrm{fine} \  | \ A_3 \bar{x}^{(i)}_\mathrm{fine} ]^\mathsf{T} \in \C^{3 \times K}$ the measurement of $x^{(i)}$ by channelwise multiplication with the system matrix $A \in \C^{3 \times K \times N_\text{fine}}$ from Section~\ref{sec:physical_model}, where $A_j \in \C^{K \times N_\mathrm{fine}}$ corresponds to the $j$-th channel of $A$,
    \item[III] $y^{\delta,(i)} = y^{(i)} + \eta^{(i)} \in \C^{3 \times K}$ the noise-perturbed measurement of $x^{(i)}$ with noise $\eta^{(i)} \in \C^{3 \times K}$, as in Section~\ref{subsec:NoiseMeas}.
\end{itemize}

For the generation of $y^{(i)}$ we select $A=A_\mathrm{opt}$, giving the best possible approximation of the measured system matrix. 

In compliance with MNIST, we propose a split of the dataset into a training set of $L_\text{train} = \num{60000}$ samples and a test set of $L_\text{test} = \num{10000}$ samples. This results in the two sets 

\begin{itemize}
    \item $\big(x^{(i)}_{\text{train}}, y^{(i)}_{\text{train}}, y^{\delta,(i)}_\text{train}\big)_{i=1,\hdots,L_\text{train}}$,
    \item $\big(x^{(i)}_{\text{test}}, y^{(i)}_{\text{test}}, y^{\delta,(i)}_\text{test}\big)_{i=1,\hdots,L_\text{test}}$.
\end{itemize}

\vspace{2.0ex}
\subsection{Reconstruction System Matrix}\label{subsubsec:rec_matrix}

To reconstruct $x^{(i)}$ from a given $y^{\delta, (i)}$, a reconstruction system matrix $A_\mathrm{rec} \in \C^{3 \times K \times N}$ needs to be selected. To prevent an inverse crime, this matrix is typically of lower spatial resolution compared to the data generation matrix $A \in \C^{3 \times K \times N_\mathrm{fine}}$, i.e., $N < N_\mathrm{fine}$. In the following we denote the spatially downsampled version of $A$ following the same physical model by $\Tilde{A}$, so that $\Tilde{A} \in \C^{3 \times K \times N}$. If the system matrix $\Tilde{A}$ differs from system matrix $A_\mathrm{rec}$ for data reconstruction, e.g. in terms of its underlying physical model, we have an operator deviation $\varepsilon: \R^{N} \longrightarrow \C^{3 \times K}$, where
\begin{equation}\label{eq:operator_dev}
    \varepsilon(x) = \Tilde{A}x - A_\mathrm{rec}x \quad \text{ for } x \in \R^{N}.
\end{equation} 
To address the operator error $\varepsilon$, oftentimes, this requires an appropriate operator correction strategy incorporated by the reconstruction method. This setting is especially relevant in the context of MPI since the system matrix of actual MPI measurements cannot be accessed. We therefore expect some kind of operator deviation when selecting a system matrix for our reconstruction methods. 

With the MPI-MNIST dataset, two different types of operator deviations are covered: A \emph{model deviation} and a \emph{noise deviation}. The former case is covered when the physical models of data generation system matrix and reconstruction system matrix differ. The noise deviation can be simulated by applying 
\begin{equation}\label{eq:H_SM}
    A_\mathrm{rec} = \Tilde{A} + H_\mathrm{SM}
\end{equation}
with some noise matrix $H_\mathrm{SM} \in \C^{3 \times K \times N}$. For all $j=1,\hdots,N$ the matrix entries $(H_\mathrm{SM})_j \in \C^{3 \times K}$ result from averaging over ten empty scanner measurements of the Bruker MPI scanner as described in Section~\ref{subsec:NoiseMeas}. Thus, by adding $H_\mathrm{SM}$, we create a pixelwise deviation of the operators given by the noise model of the MPI scanner. Both types of deviations are illustrated in Figure~\ref{fig:sms}. A combination of both deviation types will increase the operator error. Our setting aims to simulate not only the case in which a calibrated system matrix is used for reconstructing from real MPI measurements (noise deviation) but also the case where a model-based system matrix is used (model deviation). A simulation of these scenarios thus might enable the validation of operator correcting reconstruction schemes. 

\vspace{1.0ex}
\subsection{File Structure}\label{subsubsec:file_struct}

The MPI-MNIST data is provided in different files of the MDF data format~\cite{Knopp2016}, specifically tailored to MPI measurements. All data, except for the ground truth phantoms, are stored in this file format, while the ground truth phantoms are available as HDF5 files. 

\subsubsection{System Matrices}

\noindent The file \texttt{SM.tar.gz} archives twelve system matrices, covering four physical models at three different resolutions. Each matrix is stored as \texttt{/SM/SM\_\{physical model\}\_\{resolution\}.mdf}.

The physical model of the system matrix can be selected by choosing one of the contained cases for \texttt{physical model}: 
\begin{itemize}
    \item \texttt{fluid\_opt}: $A_\text{opt} \in \C^{3 \times K \times N}$ 
    \item \texttt{fluid\_small\_params}: $A_\text{SP} \in \C^{3 \times K \times N}$
    \item \texttt{fluid\_large\_params}: $A_\text{LP} \in \C^{3 \times K \times N}$
    \item \texttt{equilibrium}: $A_\text{eq} \in \C^{3 \times K \times N}$ 
\end{itemize}
The option \texttt{resolution} covers three resolutions of the spatial system matrix dimension $N$: 
\begin{itemize}
    \item \texttt{fine}: $N = N_\mathrm{fine} = 6375 \ (85  \times 75 \, \mathrm{px})$,
    \item \texttt{int}(ermediate): $N = N_\mathrm{int} = 2295 \ ( 51  \times 45 \, \mathrm{px}$),
    \item \texttt{coarse}: $N = N_\mathrm{coarse} = 255  \ ( 17  \times 15 \, \mathrm{px})$.
\end{itemize}
The \texttt{fine} resolution is used for generating the MPI measurements together with the physical model \linebreak \texttt{fluid\_opt}, while \texttt{coarse} equals the selected ground truth dimension of the MNIST phantoms. The additionally provided resolution \texttt{int} can be used in a customized setting, as discussed in Section~\ref{subsec:custom_mpi-mnist}. Furthermore, all provided system matrices can be selected for image reconstruction. In the following, we denote the system matrices by $A_p^r$ with $r = \texttt{resolution}$ and $ p = \texttt{physical model}$.

\vspace{1.0ex}
\subsubsection{Noise Measurements}\label{subsubsec:noise_meas}

\noindent The noise measurements as described in Section~\ref{subsec:NoiseMeas}  serve distinct purposes and are thus stored in seperate files. 
The data is organized as follows, where $D =\texttt{dataset}$ either takes the value \texttt{train} or \texttt{test}: 
\begin{itemize}
    \item \texttt{\{dataset\}\_noise.tar.gz} contains 
    \begin{itemize}
        \item[1.] \texttt{NoiseMeas\_phantom\_\{dataset\}.mdf} with tensor $H_{\text{phantom}}^D \in \C^{L_D \times 1 \times 3 \times K}$, $L_D$ noise samples reserved for the additive measurement noise in III,
        \item[2.] \texttt{NoiseMeas\_SM\_\{dataset\}.mdf} with tensor $H_{\text{SM}}^D \in \C^{N_\text{int}\times 1 \times 3 \times K}$, $N_\mathrm{int}$ noise samples as system matrix noise as in~\eqref{eq:H_SM},
        \item[3.] \texttt{NoiseMeas\_phantom\_bg\_\{dataset\}.mdf} with tensor $H_{\text{phantom}}^{\text{bg},D} \in \C^{100 \times 1 \times 3 \times K}$, $100$ noise samples reserved for background correction of 1,
        \item[4.] \texttt{NoiseMeas\_SM\_bg\_\{dataset\}.mdf} with tensor $H_{\text{SM}}^{\text{bg},D} \in \C^{N_\text{int}\times 1 \times 3 \times K}$, $N_\mathrm{int}$ noise samples reserved for background correction of 2.
    \end{itemize}
    \item \texttt{large\_noise.tar.gz} \\ contains \texttt{large\_NoiseMeas.mdf} with tensor $H_\text{train} \in \C^{\num{390060} \times 1 \times 3 \times K}$: This file contains \num{390060} noise samples. 
\end{itemize}

The shape of the listed tensors adapts to the conventions of the MDF file structure. The spatial dimensionality $N_\text{int}$ of 2. and 4. allows for flexibility when using the intermediate resolution. Nevertheless, only $N_\text{coarse} = 255$ spatial points are required at the coarsest resolution.

\vspace{1.0ex}
\subsubsection{Measurements}

\noindent For each part of the dataset ($D$ = \texttt{dataset = }\texttt{train} or $D$ = \texttt{dataset = }\texttt{test}), we store the MNIST phantoms and phantom measurements with the following structure:

\begin{itemize}
    \item \texttt{\{dataset\}\_obs.tar.gz} \\ contains \texttt{\{dataset\}\_obs.mdf}: Set of noise-free MPI measurements $y = [y^{(1)} | \hdots | y^{(L_D)}]\in \R^{L_D \times 1 \times 3 \times K}$ of all phantoms resulting from  $y^{(i)} = A_{\text{opt}}^{\text{fine}} \bar{x}^{(i)}_\mathrm{fine} $ as in II,
    \item \texttt{\{dataset\}\_obsnoisy.tar.gz} \\  contains \texttt{\{dataset\}\_obsnoisy.mdf}:
    Set of noise-perturbed MPI measurements $y^\delta = [y^{\delta,(1)} | \hdots | y^{\delta,(L_D)}]$ $ {\in \R^{L_D \times 1 \times 3 \times K}}$ of all phantoms resulting from $y^{\delta} = y + H_\text{phantom}^D$ as in III.
\end{itemize}

\vspace{1.0ex}
\subsubsection{Phantoms}

\noindent The phantom files are stored in $\texttt{\{dataset\}\_gt.tar.gz}$ as $\texttt{\{dataset\}\_gt.hdf5}$ for each part of the dataset (\texttt{dataset = train} or \texttt{dataset = test}). These contain the set of ground truth preprocessed and flattened MNIST phantoms $\bar{x}_\mathrm{coarse} \in \R^{L_D \times N_\mathrm{coarse}}$ with the train and test split.

\subsection{Metadata} 

The MDF file format is designed for storing MPI data together with all relevant meta information of the experimental setup. We adhered to the conventions established in~\cite{Knopp2016} for our data and made necessary adaptions.

The provided system matrices of simulated nature were adapted to a chosen calibrated system matrix. Therefore, all experimental parameters align with those of $A_\mathrm{calib}$. This also includes the signal-to-noise ratio~(SNR) stored in \texttt{/calibration/snr} of all system matrix files. 
Although the simulated system matrix does not contain noise, the SNR serves as an orientation for frequency selection during the reconstruction process, since its threshold can be seen as a regularization parameter, for further details see e.g.~\cite{Scheffler2024}. In addition to the conventional meta information, we added two parameters to the $\texttt{/tracer}$ group of \texttt{SM\_fluid\_opt\_\{resolution\}.mdf}, \texttt{SM\_fluid\_small\_params\_\{resolution\}.mdf}, \\\texttt{SM\_fluid\_large\_params\_\{resolution\}.mdf} at all three resolutions: 
\begin{itemize}
    \item \texttt{\_K\_anis}: Anisotropy constants $K_\mathrm{anis}^\mathrm{max}$ used, in $\mathrm{J}/\mathrm{m}^3$,
    \item \texttt{\_diameters}: Particle diameters $D_\mathrm{core}$ used, in $\mathrm{nm}$. 
\end{itemize}
Since \texttt{SM\_fluid\_opt\_\{resolution\}.mdf} contains multiple anisotropy constants and particle diameters, the additional parameter \texttt{\_weights} was added to the \texttt{/tracer} group and captures the weight matrix $W^\ast$ as in~\eqref{eq:lstsq} of all anisotropy constants and particle diameters. For the equilibrium model, only the parameter \texttt{\_diameters} was added to the file's meta information, since its model does not depend on any anisotropy constant.

\section{Numerical Results}\label{sec:num}

\noindent In this section we evaluate reconstructions of MPI measurements from the MPI-MNIST dataset by applying the LDA image reconstruction approach. This is complemented by a comparison of our learning-based method to standard reconstruction techniques listed in Section~\ref{subsec:standard_rec}. 

For the reconstruction of MNIST phantoms from simulated MPI measurements we choose the default MPI-MNIST setting, which selects $A = A_\mathrm{opt}^{\mathrm{fine}}$ for data generation of the measurements $(y^{(i)})_i$ together with the noise matrices $H_\mathrm{phantom}^\mathrm{train}, H_\mathrm{phantom}^\mathrm{test}$ for computing $(y^{\delta, (i)})_i$. The reconstruction of the phantoms $(x^{(i)})_i$ is performed at the coarse resolution $N_\mathrm{coarse}$. Therefore, at a concentration of $c=10\, \mathrm{mg}_\mathrm{Fe}/\mathrm{mL}$ we can simply use the data tuples $(x^{(i)}, y^{(i)}, y^{\delta, (i)})_i$ as provided on the Zenodo website, detailed in Section~\ref{sec:dataset}. For our reconstruction scenario we choose the same physical model at coarse resolution, i.e., we select $A_\mathrm{rec} = A_\mathrm{opt}^\mathrm{coarse}$ for reconstruction. For all involved components in the frequency domain, we apply a selection of frequency entries $I = \{50,\hdots,813\} \subset I_K$, where $I_L = \{0,\hdots,L-1 \}$ the index set of length $L \in \N$. This corresponds to applying a band-pass filter with a frequency range of $[77 \ \mathrm{kHz}, \num{1245} \ \mathrm{kHz}]$. Subsequently, we stack the selected frequency entries along the channel dimension, resulting in $K' = 2292$ frequency components. For $A_j \in \R^{K \times N_\mathrm{fine}}$ the $j$-th channel of $A$, we redefine $A = [A_1|_{I \times I_{N_\mathrm{fine}}}, \ A_2|_{I \times I_{N_\mathrm{fine}}}, \ A_3|_{I \times I_{N_\mathrm{fine}}}] \in \C^{K'\times N_\mathrm{fine}}$, analogously for $A_\mathrm{rec} \in C^{K' \times N_\mathrm{coarse}}$. Similarly, we define $y^{(i)} = [y_1^{(i)}|_{I}, \ y_2^{(i)}|_{I}, \ y_3^{(i)}|_{I}] \in \C^{K'}$, analogously for $y^{\delta,(i)}$, $\eta^{(i)}\in \C^{K'}$. 
For evaluation, the ground truth phantoms are retransformed to image shape, i.e., ${x^{(i)} = x_\mathrm{coarse}^{(i)} \in \R^{17 \times 15}}$. In addition to the selected concentration, we consider the same setting for concentrations $c=2,5,20,50 \,\mathrm{mg}_\mathrm{Fe}/\mathrm{mL}$. This creates five versions of the dataset. Since the noise data remains unchanged for simulating the noise-perturbed measurements $y^{\delta,(i)}$ at all concentrations, the noise level is inherently increased at a decreasing particle concentration. 

\subsection{Reconstruction Methods}

\noindent We employ four different reconstruction approaches on the noise-perturbed MPI measurements for the reconstruction of each of the MNIST phantoms. 

In addition to conventional methods introduced in Section~\ref{subsec:standard_rec}, we apply the LDA from Section~\ref{sec:lmap} by minimizing the objective in~\eqref{eq:lmap_obj} via gradient descent. In analogy to the standard methods, for the LDA we choose the regularization functional $\mathcal{R}(x) = \frac{\alpha}{2}\|x\|^2$
. Furthermore, the initial value for the minimization of~\eqref{eq:lmap_obj} is chosen as the reconstruction result obtained from the regularized Kaczmarz algorithm. Therefore, the approach can be seen as a post-processing tool. 

The regularization parameters of all reconstruction methods are selected according to an extensive grid search over the parameter spaces corresponding to the structural similarity index measure~(SSIM) as defined in~\cite{Wang2004}. The grid search is performed separately with respect to each concentration level over the first 100 images of the test dataset for the Tikhonov approach, the regularized Kaczmarz approach and the whitened regularized Kaczmarz approach for a fair comparison and representative reconstruction quality. For the LDA, we apply the grid search to the first five images of the test dataset. 

\subsection{Discrepancy Network}

 \noindent Applying the LDA requires the training of a discrepancy network. To this end, we train an invertible neural network $\varphi_\theta: \R^{2 \times K'} \longrightarrow \R^{2 \times K'}$ on noise samples $(\eta^{(i)})_i$. The network processes the real and imaginary components of the signal separately, resulting in input and output dimensions of $\R^{2 \times K'}$, where the first dimension describes the real and complex channels, respectively. The noise data from the MPI-MNIST dataset used for training originates from the samples included in the file \texttt{large\_NoiseMeas.mdf} as described in Section~\ref{sec:dataset}, where approximately $\num{300000}$ samples were used for training the network. 
 The network follows a multi-scale architecture based on the idea to capture information at different scales in order to improve the overall model's ability to represent complex data distributions. The architecture is based on affine coupling layers. 
  A detailed description of the network architecture containing approximately $13.9$ million trainable parameters can be found in Appendix~\ref{appendix:model_arch}. For approximating the true noise distribution, the loss function for the training of our INN is chosen as in~\eqref{eq:loss} and the network weights where optimized using Adam~\cite{Kingma2018} on a batch size $b= 256$ over 25 epochs of training. Note that, according to the usage of identical noise samples for the generation of all five versions of the dataset, the same trained discrepancy network can be utilized at all levels of concentration. 

\subsection{Results}

\begin{figure*}
    \centering
    \includegraphics[width=.85\textwidth]{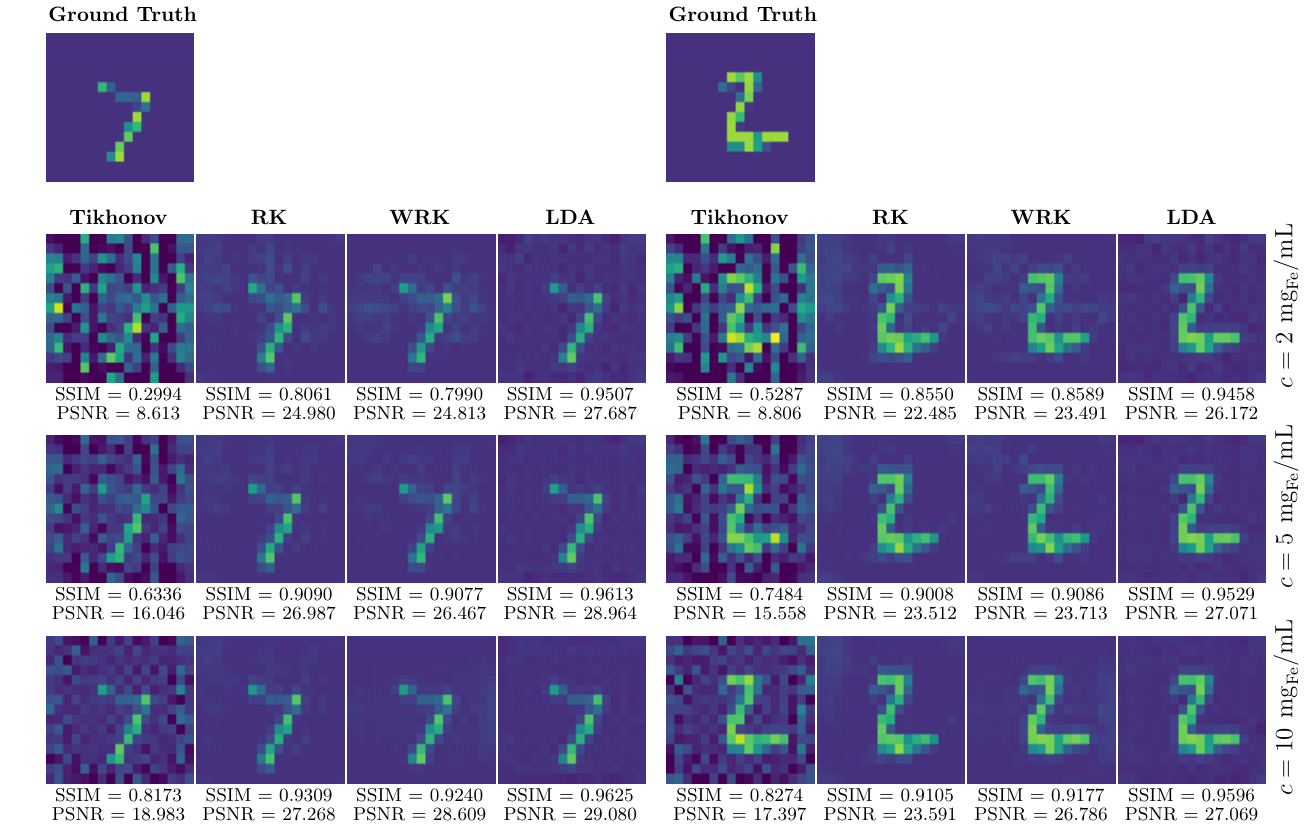}
    \caption{Reconstructions of selected phantoms of the MPI-MNIST test dataset from different reconstruction methods over concentrations $c=2,5,10 \, \mathrm{mg}_\mathrm{Fe}/\mathrm{mL}$ with corresponding error measures SSIM and PSNR compared to the ground truth.}
    \label{fig:recons_ex}
\end{figure*}

\noindent Figure~\ref{fig:recons_ex}  presents reconstructed images corresponding to selected MPI-MNIST test dataset phantoms, obtained using four distinct reconstruction methods across varying concentrations. The analytical Tikhonov approach generally yields relatively poor reconstruction quality, failing to clearly resolve the phantom structures. In contrast, the other methods demonstrate significant improvements in terms of reconstruction quality. The regularized Kaczmarz algorithm, however, introduces artifacts near the upper left corner of the images and close to the digit in the center of the image, which become increasingly visible as the concentration decreases, coinciding with higher noise levels. The whitened regularized Kaczmarz method mitigates the corner artifacts, so that they become only slightly visible, while the centered artifacts are maintained or even increased. The LDA effectively eliminates these artifacts, leaving only high-frequency background noise noticeable at concentrations $c=2,5 \, \mathrm{mg}_\mathrm{Fe}/\mathrm{mL}$. This behaviour is also reflected by both error measures, which severely differ from the standard methods especially at a low particle concentration. Among all presented methods, LDA consistently delivers superior reconstruction quality, especially at increasing noise level. This aligns with the method's design, where the noise structure learned by the discrepancy network becomes increasingly advantageous in high-noise scenarios.

\begin{figure}[t]
    \centering
    \begin{subfigure}{.45\textwidth}
        \centering
        \begin{tikzpicture}[scale=0.71]
        \begin{axis}[
            xlabel={concentration [$\mathrm{mg}_\mathrm{Fe}/\mathrm{mL}$]},
            ylabel={$\overline{\mathrm{SSIM}}$},
            legend style={at={(0.98,0.02)}, anchor=south east},
            grid=major
        ]
        
        \addplot[
            color=red,
            mark=*,
            mark options={solid}
            ]
            coordinates {
            (2, 0.4014)
            (5, 0.6911)
            (10, 0.8106)
            (20, 0.8630)
            (50, 0.8885)
            };
        
        \addlegendentry{Tikhonov}

        \addplot[
            color=blue,
            mark=*,
            mark options={solid}
            ]
            coordinates {
            (2, 0.8527)
            (5, 0.9010)
            (10, 0.9149)
            (20, 0.9365)
            (50, 0.9545)
            };
            
        \addlegendentry{RK}

        \addplot[
            color=orange,
            mark=*,
            mark options={solid}
            ]
            coordinates {
            (2, 0.8424)
            (5, 0.8872)
            (10, 0.9199)
            (20, 0.9424)
            (50, 0.9542)
            };
            
        \addlegendentry{WRK}

        \addplot[
            color=violet,
            mark=*,
            mark options={solid}
            ]
            coordinates {
            (2, 0.9314)
            (5, 0.9532)
            (10, 0.9558)
            (20, 0.9602)
            (50, 0.9709)
            };
            
        \addlegendentry{LDA}
        
        \end{axis}
        \end{tikzpicture}
        
    \subcaption{Mean SSIM values.}
    \end{subfigure}
    \hfill
    \begin{subfigure}{0.45\textwidth}
        \centering
        \begin{tikzpicture}[scale=0.71]
        \begin{axis}[
            xlabel={concentration [$\mathrm{mg}_\mathrm{Fe}/\mathrm{mL}$]},
            ylabel={$\overline{\mathrm{PSNR}}$},
            legend style={at={(0.98,0.02)}, anchor=south east},
            grid=major
        ]
        \addplot[
            color=red,
            mark=*,
            mark options={solid}
            ]
            coordinates {
            (2, 5.6092)
            (5, 13.117)
            (10, 17.944)
            (20, 21.190)
            (50, 22.975)
            };

        \addlegendentry{Tikhonov}

        \addplot[
            color=blue,
            mark=*,
            mark options={solid}
            ]
            coordinates {
            (2, 23.569)
            (5, 25.040)
            (10, 25.422)
            (20, 27.515)
            (50, 28.277)
            };
            
        \addlegendentry{RK}

        \addplot[
            color=orange,
            mark=*,
            mark options={solid}
            ]
            coordinates {
            (2, 23.641)
            (5, 24.836)
            (10, 26.995)
            (20, 27.728)
            (50, 28.041)
            };
            
        \addlegendentry{WRK}

        \addplot[
            color=violet,
            mark=*,
            mark options={solid}
            ]
            coordinates {
            (2, 25.962)
            (5, 27.335)
            (10, 27.305)
            (20, 27.981)
            (50, 28.262)
            };
            
        \addlegendentry{LDA}
        
        \end{axis}
        \end{tikzpicture}
        
    \subcaption{Mean PSNR values.}
    \end{subfigure}
    \caption{SSIM~(left) and PNSR~(right) values of the reconstruction methods evaluated over the test dataset for concentrations $c=2,5,10,20,50 \, \mathrm{mg}_\mathrm{Fe}/\mathrm{mL}$.}
    \label{fig:ssim_psnr_testset}
\end{figure}
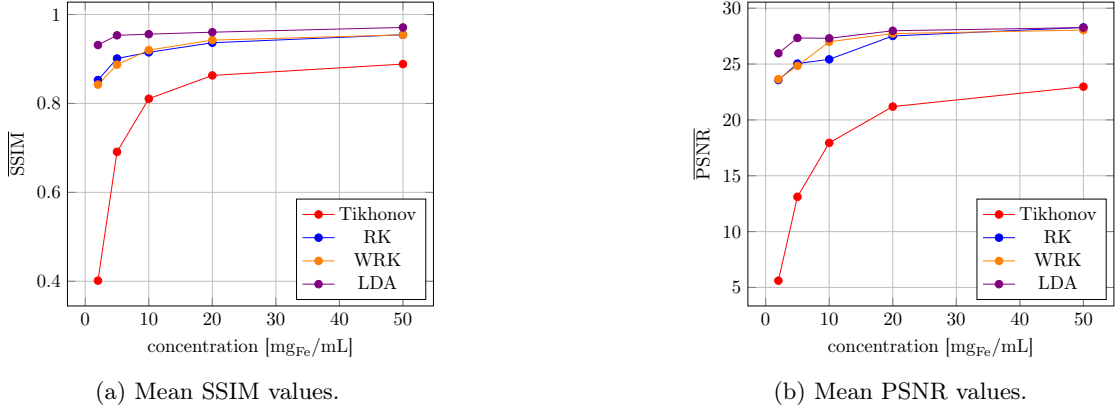

\begin{table*}[t]
\renewcommand{\arraystretch}{1.2}
\footnotesize
\centering
\begin{tabular}{cc|c|c|c|c|c|}
\cline{3-7}
                       &  & c = 2  & c = 5  & c = 10 & c = 20 & c = 50 \\ \hline
\multicolumn{1}{|c|}{\multirow{4}{*}{$\overline{\mathrm{SSIM}}$}} & Tikhonov   &  0.5087 $\pm$ 0.092 &  0.7455 $\pm$ 0.082  &  0.8000 $\pm$ 0.0818 &  0.8538 $\pm$ 0.068 & 0.8887 $\pm$ 0.048  \\ \cline{2-7} 
\multicolumn{1}{|c|}{} & RK &  0.8527 $\pm$ 0.082 & 0.9010 $\pm$ 0.050 & 0.9149 $\pm$ 0.036 & 0.9365 $\pm$ 0.039 & 0.9545 $\pm$ 0.029 \\ \cline{2-7}
\multicolumn{1}{|c|}{} & WRK & 0.8424 $\pm$ 0.082 & 0.8872 $\pm$ 0.058 & 0.9199 $\pm$ 0.051 & 0.9424 $\pm$ 0.035 & 0.9542 $\pm$ 0.026\\ \cline{2-7}
\multicolumn{1}{|c|}{} & LDA  & \textbf{0.9314 $\pm$ 0.039} & \textbf{0.9532 $\pm$ 0.027} &  \textbf{0.9558 $\pm$ 0.023} & \textbf{0.9602 $\pm$ 0.022} & \textbf{0.9709 $\pm$ 0.017} \\  \Xhline{3\arrayrulewidth}
\multicolumn{1}{|c|}{\multirow{4}{*}{$\overline{\mathrm{PSNR}}$}} & Tikhonov  & 8.633 $\pm$ 0.451&  15.587 $\pm$ 0.657 & 16.276 $\pm$ 1.553  &  19.406 $\pm$ 1.700  &  22.983 $\pm$ 1.779 \\ \cline{2-7} 
\multicolumn{1}{|c|}{} & RK  & 23.569 $\pm$ 1.129 & 25.040 $\pm$ 1.703 & 25.422 $\pm$ 1.978 &  27.525 $\pm$ 1.469 & \textbf{28.277 $\pm$ 1.459}\\ \cline{2-7}
\multicolumn{1}{|c|}{} & WRK & 23.641 $\pm$ 1.092 & 24.836 $\pm$ 1.671 & 26.995 $\pm$ 1.296 & 27.728 $\pm$ 1.527 & 28.041 $\pm$ 1.629 \\ \cline{2-7}
\multicolumn{1}{|c|}{} & LDA & \textbf{25.962 $\pm$ 1.312} & \textbf{27.335 $\pm$ 1.551} & \textbf{27.305 $\pm$ 1.666} & \textbf{27.981 $\pm$ 1.596} & 28.262 $\pm$ 1.539 \\ \hline
\end{tabular}
\caption{Average errors  over particle concentrations (in $\mathrm{mg}_\mathrm{Fe}/\mathrm{mL}$) of given reconstruction method on MNIST test set.}
\label{tab:error_meas}
\end{table*}

To investigate the general reconstruction quality on the MPI-MNIST dataset on a larger amount of phantoms, we compute the mean SSIM and PSNR values 
\begin{align*}
    \overline{\mathrm{SSIM}} &= \frac{1}{L_\mathrm{test}} \sum_{i=1}^{L_\mathrm{test}}\mathrm{SSIM}(x^{(i)}, \hat{x}^{(i)}), \\ \overline{\mathrm{PSNR}} &= \frac{1}{L_\mathrm{test}} \sum_{i=1}^{L_\mathrm{test}} \mathrm{PSNR}(x^{(i)}, \hat{x}^{(i)}),
\end{align*}
over the test dataset for different concentrations. The error measures for $c=2,5,10,20,50\, \mathrm{mg}_\mathrm{Fe}/\mathrm{mL}$ are listed in Table~\ref{tab:error_meas} and show the same effect as observed in Figure~\ref{fig:recons_ex}: The higher the noise level, the more significant the improvement of reconstruction quality of the LDA. The LDA yields the highest mean SSIM values for all concentrations, whereas the  regularized Kaczmarz algorithm performs slightly better at $c=50 \, \mathrm{mg}_\mathrm{Fe}/\mathrm{mL}$ in terms of the PSNR measure. These quantitative measures are visualized in Figure~\ref{fig:ssim_psnr_testset}, where the significance of the varying reconstruction quality over the methods at high noise levels becomes visible. All methods have the tendency to converge to a high reconstruction quality at increasing concentration, whereas the LDA shows the steepest increase towards a high reconstruction quality at low particle concentrations. However, the LDA requires a severely higher number of iterations in comparison to the other approaches, which in turn results in an increased computation time. Table~\ref{tab:iter} shows the iteration numbers of the approaches derived from our grid search together with the corresponding computation times in Table~\ref{tab:time} at $c=2,5,10,20,50\, \mathrm{mg}_\mathrm{Fe}/\mathrm{mL}$. For this study, one iteration is defined as a single gradient descent step for the LDA, or a single iteration of the regularized Kaczmarz algorithm for the RK and WRK methods. The analytical solution of the Tikhonov functional is achieved within a single reconstruction step. All computations were performed on a system with dual Intel(R) Xeon(R) Silver 4210 CPU. The computation times were averaged over the reconstruction times of the first $100$ images from the MPI-MNIST dataset. Note that we constrained our grid search for the LDA to a maximum number of $\num{2000}$ iterations, although allowing a higher iteration count for the LDA might lead to an improved reconstruction quality. To prioritize reconstruction time, the number of iterations was limited accordingly. 

\begin{table}[t]
\renewcommand{\arraystretch}{1.2}
\footnotesize
\centering
\begin{minipage}{0.45\textwidth}
    \centering
    \begin{tabular}{|c|c|c|c|c|}
    \hline
    $c\ [\mathrm{mg}_\mathrm{Fe}/\mathrm{mL}]$ & Tikhonov & RK & WRK & LDA \\ \hline
    $2$  & 1 & 1 & 1 & \num{2000} \\ \hline
    $5$  & 1 & 1 & 1 & \num{2000} \\ \hline
    $10$ & 1 & 1 & 20 & \num{2000} \\ \hline
    $20$ & 1 & 10 & 20 & \num{2000} \\ \hline
    $50$ & 1 & 200 & 20 & \num{1000} \\ \hline
    \end{tabular}
    \caption{Number of iterations for a single reconstruction of different reconstruction algorithms and concentration levels averaged over the first $100$ samples of the MPI-MNIST test dataset.}
    \label{tab:iter}
\end{minipage}
\hfill 
\begin{minipage}{0.49\textwidth}
    \centering
    \begin{tabular}{|c|c|c|c|c|}
    \hline
    $c\ [\mathrm{mg}_\mathrm{Fe}/\mathrm{mL}]$ & Tikhonov & RK & WRK & LDA \\ \hline
    $2$  & 0.0079 & 0.0665 & 0.0691 & 89.5331 \\ \hline
    $5$  & 0.0078 & 0.0684 & 0.0701 & 82.9576 \\ \hline
    $10$ & 0.0081 & 0.0512 & 0.9854 & 84.8896 \\ \hline
    $20$ & 0.0080 & 0.4830 & 0.9655 & 87.0539 \\ \hline
    $50$ & 0.0070 & 10.0687 & 1.0601 & 42.2373 \\ \hline
    \end{tabular}
    \caption{Time in s for a single reconstruction of different reconstruction algorithms and concentration levels averaged over the first $100$ samples of the MPI-MNIST test dataset.}
    \label{tab:time}
\end{minipage}
\end{table}

\section{Conclusion}

\noindent In the present work, we addressed two challenges in MPI image reconstruction: The formulation of a reconstruction algorithm that accurately captures the structure of MPI scanner noise and the scarcity of available MPI datasets.

From an algorithmic perspective, we introduced the LDA, which leverages the features and capabilities of the MPI-MNIST dataset. Unlike conventional Tikhonov regularization, our approach relaxes the Gaussian noise assumption
. The applied MAP estimation within the LDA includes a formulation of the likelihood of the data under the model which comprises a learned version of the true noise distribution by an invertible neural network. Therefore, the method is capable of taking complex noise structures into account. This incorporation of problem-specific noise features directly addresses the noise aspect of the MPI modeling chain.

Addressing the second challenge, we presented the publicly available MPI-MNIST dataset, which is beneficial for the development and testing of reconstruction methods. 
The extensive collection of ground truth phantoms and corresponding simulated MPI measurements as well as a substantial number of additional noise samples captured by the preclinical MPI scanner allow for a detailed analysis of complex data structures. The available data scope facilitates the application of learning-based approaches. Moreover, the flexibility of the dataset enables the generation of customized data versions and a selection of various reconstruction scenarios. 
 
In our numerical experiments we demonstrated the effectiveness of our reconstruction approach by training the discrepancy network on the noise data provided with the MPI-MNIST dataset. Compared to classical reconstruction methods, 
our learned approach showed significant improvements in reconstruction quality, particularly at higher noise levels, which occurs with a reduction of the particle concentration. 
However, a current limitation is its extended reconstruction duration in comparison to classical methods.

These findings serve as a basis for several future directions. Firstly, it is essential to evaluate the performance of the LDA on real MPI measurements by using the trained model from our numerical experiments. This step is crucial for validating the method’s generalizability in practical scenarios. Additionally, the regularization properties of the LDA remain unexplored. 
Examining the conditions of the discrepancy network, including specific network architecture choices, under which the approach maintains existence, stability, and convergence of a regularization scheme will be an essential focus for future research. In this work, we focus on the discrepancy functional within the Tikhonov objective by addressing the noise distribution. However, the approach could be extended by additionally relaxing the Gaussian assumption on the prior. In such cases, we could introduce an additional network to learn the prior distribution, potentially leading to a more accurate penalty term and bringing the method even closer to a realistic MPI setting.

\section*{Acknowledgements}
\noindent The authors would like to thank Florian Thieben and Konrad Scheffler for conducting the experiments. Furthermore, the authors would like to thank Marija Boberg for her useful assistance with data preparation.
Meira Iske acknowledges the support of the Deutsche Forschungsgemeinschaft
(DFG, German Research Foundation) - Project number 281474342/GRK2224/2. Hannes Albers acknowledges financial support from the KIWi project funded by the German Federal Ministry of Education and Research (BMBF) - Project number 05M22LBA. 

\printbibliography

@misc{lecun1998mnist,
  title={The {MNIST} database of handwritten digits},
  author={LeCun, Y},
  url={http://yann.lecun.com/exdb/mnist/},
  year={1998}
}

@ARTICLE{maass2024,
  author={Maass, M and Kluth, T and Droigk, C and Albers, H and Scheffler, K and Mertins, A and Knopp, T},
  journal={IEEE Transactions on Computational Imaging}, 
  title={Equilibrium model with anisotropy for model-based reconstruction in magnetic particle imaging}, 
  year={2024},
  volume={10},
  number={},
  pages={1588-1601},
  doi={10.1109/TCI.2024.3490381}}

@article{gungor2023deq,
  title={{DEQ-MPI}: A deep equilibrium reconstruction with learned consistency for magnetic particle imaging},
  author={G{\"u}ng{\"o}r, A and Askin, B and Soydan, D A and Top, C B and Saritas, E U and {\c{C}}ukur, T},
  journal={IEEE Transactions on Medical Imaging},
  volume={43},
  number={1},
  pages={321--334},
  year={2023},
  publisher={IEEE}
}

@misc{Knopp2016,
  author = {Knopp, T and Viereck, T and Bringout, G and Ahlborg, M and Rahmer, J and Hofmann, M},
  title = {{MDF}: Magnetic particle imaging data format},
  doi = {10.48550/arxiv.1602.06072},
  url = {https://arxiv.org/abs/1602.06072},
  publisher = {arXiv},
  year = {2020}}

@article{Kluth2019,
doi = {10.1088/1361-6560/ab1a4f},
year = {2019},
publisher = {IOP Publishing},
volume = {64},
number = {12},
pages = {125026},
author = {T Kluth and B Jin},
title = {Enhanced reconstruction in magnetic particle imaging by whitening and randomized {SVD} approximation},
journal = {Physics in Medicine \& Biology}
}

@article{Kluth2019_towards,
doi = {10.1088/1367-2630/ab4938},
year = {2019},
publisher = {IOP Publishing},
volume = {21},
number = {10},
pages = {103032},
author = {T Kluth and P Szwargulski and T Knopp},
title = {Towards accurate modeling of the multidimensional magnetic particle imaging physics},
journal = {New Journal of Physics}
}

@article{Dax1993,
author = {Dax, A},
title = {On row relaxation methods for large constrained least squares problems},
journal = {SIAM Journal on Scientific Computing},
volume = {14},
number = {3},
pages = {570-584},
year = {1993},
doi = {10.1137/0914036}
}

@inproceedings{Kingma2018,
 author = {Kingma, D and Dhariwal, P},
 booktitle = {Advances in neural information processing systems},
 editor = {S. Bengio and H. Wallach and H. Larochelle and K. Grauman and N. Cesa-Bianchi and R. Garnett},
 pages = {},
 publisher = {Curran Associates, Inc.},
 title = {Glow: Generative flow with invertible 1x1 convolutions},
 volume = {31},
 year = {2018}
}

@article{Wang2004,
  author={Z Wang and Bovik, A and Sheikh, H and Simoncelli, E},
  journal={IEEE Transactions on Image Processing}, 
  title={Image quality assessment: From error visibility to structural similarity}, 
  year={2004},
  volume={13},
  number={4},
  pages={600-612},
  doi={10.1109/TIP.2003.819861}
}

@article{Gleich2005,
 author={B Gleich and J Weizenecker},
 title={Tomographic imaging using the nonlinear response of magnetic particles},
 journal={Nature},
 year={2005},
 doi={10.1038/nature03808},
 volume={435},
 pages={1214-7},
}

@article{Knopp2017,
doi = {10.1088/1361-6560/aa6c99},
year = {2017},
publisher = {IOP Publishing},
volume = {62},
number = {14},
pages = {R124},
author = {T Knopp and N Gdaniec and M Möddel},
title = {Magnetic particle imaging: From proof of principle to preclinical applications},
journal = {Physics in Medicine \& Biology}
}

@article{Weizenecker2009,
doi = {10.1088/0031-9155/54/5/L01},
year = {2009},
publisher = {IOP Publishing},
volume = {54},
number = {5},
pages = {L1},
author = {J Weizenecker and B Gleich and J Rahmer and H Dahnke and J Borgert},
title = {Three-dimensional real-time in vivo magnetic particle imaging},
journal = {Physics in Medicine \& Biology}
}

@article{Haegele2012,
author = {Haegele, J and Rahmer, J and Gleich, B and Borgert, J and Wojtczyk, H and Panagiotopoulos, N and Buzug, T and Barkhausen, J and Vogt, F.},
title = {Magnetic particle imaging: Visualization of instruments for cardiovascular intervention},
journal = {Radiology},
volume = {265},
number = {3},
pages = {933-938},
year = {2012},
doi = {10.1148/radiol.12120424}
}

@article{Herz2019,
author = {S Herz and P Vogel and T Kampf and P Dietrich and S Veldhoen and M A. Rückert and R Kickuth and V Behr and T Bley},
title ={Magnetic particle imaging–guided stenting},
journal = {Journal of Endovascular Therapy},
volume = {26},
number = {4},
pages = {512-519},
year = {2019},
doi = {10.1177/1526602819851202}
}

@article{Gleich2008,
doi = {10.1088/0031-9155/53/6/N01},
year = {2008},
publisher = {IOP Publishing},
volume = {53},
number = {6},
pages = {N81},
author = {B Gleich and J Weizenecker and J Borgert},
title = {Experimental results on fast {2D}-encoded magnetic particle imaging},
journal = {Physics in Medicine \& Biology}
}

@article{Yu20217,
author = {Yu, E and Bishop, M and Zheng, B and Ferguson, R and Khandhar, A and Kemp, S and Krishnan, K and Goodwill, P and Conolly, S},
title = {Magnetic particle imaging: A novel in vivo imaging platform for cancer detection},
journal = {Nano Letters},
volume = {17},
number = {3},
pages = {1648-1654},
year = {2017},
doi = {10.1021/acs.nanolett.6b04865}
}

@article{Sehl2020,
    author = {O Sehl and J Gevaert and K Meli and N Knie and P Foster}, 
    title = {A perspective on cell tracking with magnetic particle imaging},
    journal = {Tomography},
    year = {2020},
    pages={315-324},
    volume={6},
    number={4}, 
    doi={10.18383/j.tom.2020.00043}
}

@article{Bulte2015,
    author = {J Bulte and P Walczak and M Janowski and K Krishnan and H Arami and A Halkola and B Gleich and J Rahmer},
    title = {Quantitative "hot spot" imaging of transplanted stem cells using superparamagnetic tracers and magnetic particle imaging ({MPI})},
    journal ={Tomography},
    year = {2015},
    volume={1},
    number={2},
    pages={91-97},
    doi={10.18383/j.tom.2015.00172}
}

@article{Herz2018,
  author={Herz, S and Vogel, P and Kampf, T and Rückert, M and Veldhoen, S and Behr, V and Bley, T },
  journal={IEEE Transactions on Medical Imaging}, 
  title={Magnetic particle imaging for quantification of vascular stenoses: A phantom study}, 
  year={2018},
  volume={37},
  number={1},
  pages={61-67},
  doi={10.1109/TMI.2017.2717958}
}

@article{Rahmer2017,
  author={Rahmer, J and Wirtz, D and Bontus, C and Borgert, J and Gleich, B},
  journal={IEEE Transactions on Medical Imaging}, 
  title={Interactive magnetic catheter steering with {3-D} real-time feedback using multi-color magnetic particle imaging}, 
  year={2017},
  volume={36},
  number={7},
  pages={1449-1456},
  doi={10.1109/TMI.2017.2679099}
}

@article{Buzug2012,
title = {Magnetic particle imaging: Introduction to imaging and hardware realization},
journal = {Zeitschrift für Medizinische Physik},
volume = {22},
number = {4},
pages = {323-334},
year = {2012},
doi = {10.1016/j.zemedi.2012.07.004},
author = {T Buzug and G Bringout and M Erbe and K Gräfe and M Graeser and M Grüttner and A Halkola and T Sattel and W Tenner and H Wojtczyk and J Haegele and F Vogt and J Barkhausen and K Lüdtke-Buzug}
}

@article{Knopp2020,
title = {{OpenMPIData}: An initiative for freely accessible magnetic particle imaging data},
journal = {Data in Brief},
volume = {28},
pages = {104971},
year = {2020},
doi = {10.1016/j.dib.2019.104971},
author = {T Knopp and P Szwargulski and F Griese and M Gräser},
}

@article{Knopp2010,
doi = {10.1088/0031-9155/55/6/003},
year = {2010},
publisher = {IOP Publishing},
volume = {55},
number = {6},
pages = {1577},
author = {T Knopp and J Rahmer and T  Sattel and S Biederer and J Weizenecker and B Gleich and J Borgert and T  Buzug},
title = {Weighted iterative reconstruction for magnetic particle imaging},
journal = {Physics in Medicine \& Biology}
}

@article{Knopp2010_model,
  author={Knopp, T and Sattel, T and Biederer, S and Rahmer, J and Weizenecker, J and Gleich, B and Borgert, J and Buzug, T},
  journal={IEEE Transactions on Medical Imaging}, 
  title={Model-based reconstruction for magnetic particle imaging}, 
  year={2010},
  volume={29},
  number={1},
  pages={12-18},
  doi={10.1109/TMI.2009.2021612}}

@article{Knopp2010_model2D,
author = {Knopp, T and Biederer, S and Sattel, T and Rahmer, J and Weizenecker, J and Gleich, B and Borgert, J and Buzug, T},
title = {{2D} model-based reconstruction for magnetic particle imaging},
journal = {Medical Physics},
volume = {37},
number = {2},
pages = {485-491},
doi = {10.1118/1.3271258},
year = {2010}
}

@article{Dittmer2021,
    author ={S Dittmer and T Kluth and M Henriksen and P Maass},
    title ={Deep image prior for {3D} magnetic particle imaging: A quantitative comparison of regularization techniques on Open {MPI} dataset},
    journal ={International Journal on Magnetic Particle Imaging},
    year ={2021}, 
    volume={7},
    number={1},
    doi={10.18416/IJMPI.2021.2103001 }
}

@article{Kluth2020,
    author = {T Kluth},
    title = {L1 data fitting for robust reconstruction in magnetic particle imaging: Quantitative evaluation on {Open MPI} dataset},
    journal = {International Journal on Magnetic particle Imaging},
    year = {2020},
    volume={6},
    number={2},
    doi={10.18416/IJMPI.2020.2012001}
}

@article{Shang2022,
doi = {10.1088/1361-6560/ac6e24},
year = {2022},
publisher = {IOP Publishing},
volume = {67},
number = {12},
pages = {125012},
author = {Y Shang and J Liu and L Zhang and X Wu and P Zhang and L Yin and H Hui and J Tian},
title = {Deep learning for improving the spatial resolution of magnetic particle imaging},
journal = {Physics in Medicine \& Biology}
}

@inproceedings{Knopp2008,
  author={Knopp, T and Biederer, S and Sattel, T and Buzug, T},
  booktitle={2008 IEEE Nuclear Science Symposium Conference Record}, 
  title={Singular value analysis for magnetic particle imaging}, 
  year={2008},
  volume={},
  number={},
  pages={4525-4529},
  doi={10.1109/NSSMIC.2008.4774296}}

@article{Knopp2016_2,
doi = {10.1088/0031-9155/61/11/N257},
year = {2016},
publisher = {IOP Publishing},
volume = {61},
number = {11},
pages = {N257},
author = {T Knopp and M Hofmann},
title = {Online reconstruction of {3D} magnetic particle imaging data},
journal = {Physics in Medicine \& Biology}
}

@article{Kluth2017,
    author = {T Kluth and P Maass},
    title = {Model uncertainty in magnetic particle imaging: Nonlinear problem formulation and model-based sparse reconstruction},
    journal = {International Journal of Magnetic Particle Imaging},
    year =  {2017},
    doi = {10.18416/ijmpi.2017.1707004},
    volume={3},
    number={2}
}

@article{Blanke2020,
doi = {10.1088/1361-6420/abb5e1},
year = {2020},
publisher = {IOP Publishing},
volume = {36},
number = {12},
pages = {124001},
author = {S Blanke and B Hahn and A Wald},
title = {Inverse problems with inexact forward operator: Iterative regularization and application in dynamic imaging},
journal = {Inverse Problems}
}

@article{Nitzsche2022,
    author = {M Nitzsche and H Albers and B Hahn and T Kluth},
    title = {Compensating model imperfections during
                image reconstruction via {RESESOP}},
    journal = {International Journal of Magnetic Particle Imaging},
    year =  {2022},
    doi = {10.18416/ijmpi.2022.2203062},
    volume={8},
    number={1}
}

@article{Bathke2017,
    author = {C Bathke and T Kluth and C Brandt and P Maass},
    title = {Improved image reconstruction in magnetic particle imaging using structural a priori information},
    journal = {International Journal of Magnetic Particle Imaging},
    year =  {2017},
    doi = {10.18416/ijmpi.2017.1703015},
    volume={3},
    number={1}
}

@article{Storath2017,
    author = {M Storath and C Brandt and M Hofmann and T Knopp and J Salamon and A Weber and A Weinmann},
    title = {Edge preserving and noise reducing reconstruction for magnetic particle imaging},
    journal = {IEEE Trans Med Imaging},
    year = {2017},
    volume={36},
    number={1},
    doi={10.1109/TMI.2016.2593954}
}

@article{Thieben2024transfer, 
    title={Experimental parameter calibration of the scanner model for model-based {MPI}}, 
    doi={10.18416/IJMPI.2024.2403025}, 
    journal={International Journal on Magnetic Particle Imaging IJMPI}, author={Thieben, F and Albers, H and Mohn, F and Foerger, F and Boberg, M and Scheffler, K and Möddel, M and Kluth, T and Knopp, T}, 
    year={2024}
}

@article{Knopp2023,
    author = {T Knopp and P Jürß and M Grosser},
    title = {A deep learning approach for automatic image reconstruction in {MPI}},
    journal = {International Journal of Magnetic Particle Imaging},
    year =  {2023},
    doi = {10.18416/ijmpi.2023.2303008},
    volume={9},
    number={1}
}

@article{Gladiss2022,
    author = {A von Gladiss and I Kramer and N Theisen and R Memmesheimer and A Bakenecker and T Buzug  and D Paulus},
    title = {Data augmentation for training a neural network for image reconstruction in {MPI}},
    journal = {International Journal of Magnetic Particle Imaging},
    year =  {2022},
    doi = {10.18416/ijmpi.2022.2203058},
    volume={8},
    number={1}
}

@inproceedings{Gungor2024,
    author = {A Güngör and E Saritas and T Çukur},
    title = {A deep equilibrium technique for {3D MPI} reconstruction},
    booktitle = {International Journal of Magnetic Particle Imaging},
    year = {2024},
    doi = {10.18416/ijmpi.2024.2403009},
    volume={10},
    number={1}
}

@InProceedings{Askin2022,
author="Askin, B
and G{\"u}ng{\"o}r, A
and Alptekin Soydan, D
and Saritas, E
and Top, C
and Cukur, T",
title="{PP-MPI}: A deep plug-and-play prior for magnetic particle imaging reconstruction",
booktitle="Machine Learning for Medical Image Reconstruction",
year="2022",
publisher="Springer International Publishing",
pages="105-114"
}

@inproceedings{Gungor2023,
    author = {A Güngör and B Askin and D Soydan and E Saritas and C Top and T Çukur},
    title = {A denoiser scaling technique for plug-and-play {MPI} reconstruction},
    booktitle = {International Journal of Magnetic Particle Imaging},
    year = {2023},
    doi = {10.18416/ijmpi.2023.2303041},
    volume = {9},
    numer={1}
}

@inproceedings{Tsanda2024,
    author = {A Tsanda and P Jürß and N Hackelberg and M Grosser and M Möddel and T Knopp},
    title = {Extension of the Kaczmarz algorithm with a deep plug-and-play regularizer},
    booktitle = {International Journal of Magnetic Particle Imaging},
    year = {2024},
    doi = {10.18416/ijmpi.2024.2403010},
    volume = {10},
    numer={1}
}

@InProceedings{Gladiss2022_1,
author="von Gladiss, A
and Memmesheimer, R
and Theisen, N
and Bakenecker, A
and Buzug, T
and Paulus, D",
editor="Maier-Hein, Klaus
and Deserno, Thomas M.
and Handels, Heinz
and Maier, Andreas
and Palm, Christoph
and Tolxdorff, Thomas",
title="Reconstruction of {1D} images with a neural network for magnetic particle imaging",
booktitle="Bildverarbeitung f{\"u}r die Medizin 2022",
year="2022",
publisher="Springer Fachmedien Wiesbaden",
pages="247-252",
doi="10.1007/978-3-658-36932-3_52"
}

@inproceedings{Koch2019,
    author = {P Koch and M Maass and M Bruhns and C Droigk and T Parbs and
A Mertins},
    title = {Neural network for reconstruction of {MPI} images},
    booktitle = {9th International Workshop on Magnetic Particle Imaging (IWMPI)},
    year = {2019} 
}

@inproceedings{Dinh2017,
    author = {L Dinh and J Sohl-Dickstein and S Bengio},
    title = {Density estimation using Real {NVP}},
    booktitle = {ICLR},
    year = {2017}
}

@article{Scheffler2024,
doi = {10.1088/1361-6560/ad2231},
year = {2024},
publisher = {IOP Publishing},
volume = {69},
number = {4},
pages = {045024},
author = {K Scheffler and M Boberg and T Knopp},
title = {Solving the {MPI} reconstruction problem with automatically tuned regularization parameters},
journal = {Physics in Medicine \& Biology}
}

@book{Kaipio2004,
    author = {J Kaipio and E Somersalo},
    title = {Statistical and Computational Inverse Problems},
    publisher ={Springer New York, NY},
    year = {2004}
}

@article{Albers2022,
title = {Modeling the magnetization dynamics for large ensembles of immobilized magnetic nanoparticles in multi-dimensional magnetic particle imaging},
journal = {Journal of Magnetism and Magnetic Materials},
volume = {543},
pages = {168534},
year = {2022},
doi = {10.1016/j.jmmm.2021.168534},
author = {H Albers and T Knopp and M Möddel and M Boberg and T Kluth}
}

@article{Albers2022simulating,
title = {Simulating magnetization dynamics of large ensembles of single domain nanoparticles: Numerical study of Brown/Néel dynamics and parameter identification problems in magnetic particle imaging},
journal = {Journal of Magnetism and Magnetic Materials},
volume = {541},
pages = {168508},
year = {2022},
doi = {10.1016/j.jmmm.2021.168508},
author = {H Albers and T Kluth and T Knopp}
}

@InProceedings{Yoshida2012,
author="Yoshida, T
and Enpuku, K",
editor="Buzug, Thorsten M.
and Borgert, J{\"o}rn",
title="Nonlinear behavior of magnetic fluid in Brownian relaxation: Numerical simulation and derivation of empirical model",
booktitle="Magnetic Particle Imaging",
year="2012",
publisher="Springer Berlin Heidelberg"
}

@article{Kluth2018,
doi = {10.1088/1361-6420/aac535},
year = {2018},
publisher = {IOP Publishing},
volume = {34},
number = {8},
pages = {083001},
author = {T Kluth},
title = {Mathematical models for magnetic particle imaging},
journal = {Inverse Problems}
}

@article{Them2015,
    author = {Them, K and Kaul, M and Jung, C and Hofmann, M and Mummert, T and Werner, F and Knopp, T} ,
    title = {Sensitivity enhancement in magnetic particle imaging by background subtraction} ,
    journal ={IEEE Transactions on Medical Imaging},
    year = {2015},
    volume={35},
    number ={3},
    pages = {893--900},
    doi={10.1109/TMI.2015.2501462}
}

@article{Knopp2019_bg,
doi = {10.1088/1361-6560/ab2480},
year = {2019},
publisher = {IOP Publishing},
volume = {64},
number = {12},
pages = {125013},
author = {T Knopp and N Gdaniec and R Rehr and M Graeser and T Gerkmann},
title = {Correction of linear system drifts in magnetic particle imaging},
journal = {Physics in Medicine \& Biology}
}

@article{Pinetz2021,
author = {Pinetz, T and Kobler, E and Pock, T and Effland, A},
title = {Shared prior learning of energy-based models for image reconstruction},
journal = {SIAM Journal on Imaging Sciences},
volume = {14},
number = {4},
pages = {1706--1748},
year = {2021},
doi = {10.1137/20M1380016}
}

@inbook{Calatroni2017,
title = {8. Bilevel approaches for learning of variational imaging models},
booktitle = {Variational Methods: In Imaging and Geometric Control},
author = {L Calatroni and C Cao and J C De los Reyes and C-B Schönlieb and T Valkonen},
editor = {M Bergounioux and G Peyré and C Schnörr and J-B Caillau and T Haberkorn},
publisher = {De Gruyter},
address = {Berlin, Boston},
pages = {252--290},
doi = {10.1515/9783110430394-008},
year = {2017}
}

@article{bulte2022vivo,
  title={In vivo cellular magnetic imaging: Labeled versus unlabeled cells},
  author={Bulte, Jeff WM and Wang, Chao and Shakeri-Zadeh, Ali},
  journal={Advanced functional materials},
  volume={32},
  number={50},
  pages={2207626},
  year={2022},
  publisher={Wiley Online Library}
}

@article{kaul2018magnetic,
  title={Magnetic particle imaging for in vivo blood flow velocity measurements in mice},
  author={Kaul, Michael G and Salamon, Johannes and Knopp, Tobias and Ittrich, Harald and Adam, Gerhard and Weller, Horst and Jung, Caroline},
  journal={Physics in Medicine \& Biology},
  volume={63},
  number={6},
  pages={064001},
  year={2018},
  publisher={IOP Publishing}
}

@article{bauer2016high,
  title={High-performance iron oxide nanoparticles for magnetic particle imaging--guided hyperthermia {(hMPI)}},
  author={Bauer, L M and Situ, S F and Griswold, M A and Samia, A C S},
  journal={Nanoscale},
  volume={8},
  number={24},
  pages={12162--12169},
  year={2016},
  publisher={Royal Society of Chemistry}
}

@article{ahlborg2022first,
  title={First dedicated balloon catheter for magnetic particle imaging},
  author={Ahlborg, M and Friedrich, T and G{\"o}ttsche, T and Scheitenberger, V and Linemann, R and Wattenberg, M and Buessen, A T and Knopp, T and Szwargulski, P and Kaul, M G and others},
  journal={IEEE transactions on medical imaging},
  volume={41},
  number={11},
  pages={3301--3308},
  year={2022},
  publisher={IEEE}
}

\appendix

\section{Modifications of MPI-MNIST}\label{subsec:custom_mpi-mnist}

\noindent The fully provided dataset described in Section~\ref{sec:dataset} contains preselected parameters such as a fixed particle concentration, reconstruction resolution, as well as a data-generating system matrix. However, together with the provided data, this framework can be customized. 

\subsection*{Particle Concentration}

\noindent For changing the particle concentration to a concentration $\Tilde{c}$, all steps of the preprocessing in Section~\ref{subsubsec:prep_MNIST} can be applied to the phantoms, whereas the pixel values are scaled by the desired concentration $\Tilde{c}$, so that $x_{j,k} \in [0,\Tilde{c}]$ for all $j=1,\hdots,W_\mathrm{coarse}$ and $k=1,\hdots,H_\mathrm{coarse}$. The measurements $y$ can either be scaled by the factor $\nicefrac{\Tilde{c}}{c}$ or $y$ and $y^\delta$ can be computed as in I. and II., respectively. By changing the concentration, the noise level of $y^\delta$ will inherently change -- a higher tracer concentration will lead to a larger signal and thus a higher SNR when the noise is kept fixed. 

\subsection*{Resolution and physical model}

\noindent The resolution and underlying physical model can be adapted both for the data generation and data reconstruction. In both cases, the system matrix is selected as 
 $A_p^r$ for $r=\texttt{coarse},$ $ \texttt{int},$ $ \texttt{fine}$ and $p = \texttt{fluid\_opt}$,  $\texttt{fluid\_small\_params},$ $\texttt{fluid\_large\_params},$ $\texttt{equilibrium}$ from the provided system matrix files listed in the previous section. In order to obtain either the corresponding ground truth or data generation phantoms $x_r, \bar{x}_r$ of the desired resolution $r$, 
 the following steps are required:
 \begin{itemize}
     \item[(i)]  $r=\texttt{fine}$: Apply preprocessing as in Section~\ref{subsubsec:prep_MNIST} to obtain $x^{(i)}_\text{fine}$ and $\bar{x}^{(i)}_\text{fine}$. 
     \item[(ii)]  $r=\texttt{int}$: Change the upsampling factor $s$ in Section~\ref{subsubsec:prep_MNIST} to $s = 3$ so that $x^{(i)}_\mathrm{int} \in \R^{51 \times 45}$. Subsequently flatten analogously to obtain $\bar{x}^{(i)}_\mathrm{int} \in \R^{N_\mathrm{int}}$.
     \item[(iii)] $r=\texttt{coarse}$: Set $\bar{x}^{(i)}_r = \bar{x}^{(i)}_\mathrm{coarse}$ and  obtain $x^{(i)}_\mathrm{coarse} \in \R^{17 \times 15}$ by reversing the matrix to vector transformation described in Section~\ref{subsubsec:prep_MNIST}.
 \end{itemize}
For the data generation process, compute $y^{(i)} = A_p^r \bar{x}^{(i)}_r$ and $y^{\delta,(i)} = y^{(i)} + \eta^{(i)}$, where $\eta^{(i)}$ results from $H_\mathrm{phantom}^D$.

\section{Model Architecture}\label{appendix:model_arch}

\noindent The model architecture used in Section~\ref{sec:num} follows a multi-scale architecture inspired by~\cite{Dinh2017} that comprises affine coupling layers with shape transforming layers. 
The real part $\eta^\mathrm{real}$ and imaginary part $\eta^\mathrm{imag}$ of our network input were considered channelwise, so that $\eta = [\eta^\mathrm{real} \ | \ \eta^\mathrm{imag}] \in \R^{p \times K'}$ with $p=2$ channels. We define by $\eta_j := [\eta^\mathrm{real}_j \ | \ \eta^\mathrm{imag}_j]$ the $j$-th component of each channel with $j=0,\hdots,K'-1$. By treating each channel separately, an affine coupling layer $f_\mathrm{ac}:\R^{p \times K'} \longrightarrow \R^{p \times K'}$ 
is defined by  

\small
\begin{equation*}
    \left(f_\mathrm{ac}(\eta)\right)_j = 
    \begin{cases} 
        \eta_j, & \text{if } j \in I_{d,K'}, \\
        \eta_j \cdot \exp\big(s((\eta_j)_{j \in I_{d,K'}})\big) + t\big((\eta_j)_{j \in I_{d,K'}}\big), & \text{if } j \in I_{d,K'}^C,
    \end{cases}
\end{equation*}
\normalsize
for $0 < d \leq K'-1$ and an index set $I_{d,K'} \subseteq \{ j \ : \ 0 \leq j \leq K'-1 \}$ with $| I_{d,K'} | = d$ and $I_{d,K'}^C = I_{K'} \setminus I_{d,K'}$. The functions $s,t:\R^{d} \longrightarrow \R^d$ are neural networks, referred to as \emph{scaling} and \emph{shifting networks}. The inverse of $y = f_\mathrm{ac}(\eta)$ can be computed by 

\small
\begin{equation*}
    \eta_j = 
    \begin{cases} 
        y_j, & \text{if } j \in I_{d,K'}, \\[0.5mm]
        \big(y_j - t\big((y_j)_{j \in I_{d,K'}}\big)\big) \cdot \exp\big(-s\big((y_j)_{j \in I_{d,K'}}\big)\big),
        & \text{if } j \in I_{d,K'}^C,
    \end{cases}
\end{equation*}
\normalsize
without the necessity to compute the inverse of $s,t$. Its Jacobian is a lower triangular matrix, 
with a tractable Jacobian determinant. To operate on different scales, we incorporate squeeze layers, split layers and concatenation layers in addition to the coupling layers. These are responsible for transforming the data shapes. The squeeze layers $f_\mathrm{squeeze}: \R^{p \times K'} \longrightarrow \R^{2p \times K'/2}$ compress the frequency dimensions while expanding the number of channels, transforming the input $w$ by
\begin{equation*}
    \left(f_\mathrm{squeeze}(w)\right)_{j,k} = w_{\lceil j/p \rceil, \ k \text{ mod } p}, 
\end{equation*}
 for $j = 0,\hdots,2p-1$, $k = 0,\hdots,K'/2-1$, as illustrated in Figure~\ref{fig:squeeze_split}. Split layers $f_\mathrm{split}: \R^{p\times K'} \longrightarrow \R^{p/2 \times K'} \times \R^{p/2 \times K'}$ with ${f_\mathrm{split}(w) = (w_{0:p/2-1, 0:K'-1}, w_{p/2:p-1, 0:K'-1})}$, on the other hand, partition an input along its channel dimension. This enables separate processing pathways that are later recombined to restore the original dimensionality~(see Figure~\ref{fig:squeeze_split}). The concatenation layer ${f_\mathrm{concat}: \R^{p/2 \times K'} \times \R^{p/2 \times K} \longrightarrow \R^{p \times K}}$, used as an output layer of the network, is responsible for this recombination by concatenating two inputs along the channel dimension, where
\begin{equation*}
    \left(f_\mathrm{concat}(v,w) \right)_{j,k} = \begin{cases}
        v_{j,k} \quad &\text{if } j\leq p/2-1 \\
        w_{j,k} &\text{otherwise}.
    \end{cases}
\end{equation*} 
Figure~\ref{fig:network_architecture} contains the detailed overall layer architecture, which extends to a third dimension referring to the batch size $b$. 
Table~\ref{tab:st_arch} shows the chosen identical architectures of the shifting and the scaling networks within each affine coupling layer $\ell = 1,\hdots,8$. The width of the fully connected layers depends on the scale of the coupling block with 
\begin{equation*}
    W_\ell = \begin{cases}
        512 & \text{ for } \ell = 1,2 \\
        256 & \text{ for } \ell = 3,4 \\
        128 & \text{ for } \ell = 5,\hdots,8. \\
    \end{cases}
\end{equation*}
This network design allows for efficient transformation of the input across multiple scales, with the affine coupling layers ensuring invertibility and computationally efficiency of the network. For further reading on the layer types we refer to~\cite{Dinh2017}.

In our setting, the index set $I_{d,K'}$ 
was implemented by defining $I_{d,K'} := \{j \ : \ j \mod 2 = 0, \ 0 \leq j \leq K'-1 \}$ with $d=\lfloor \nicefrac{K'}{2}\rfloor$, therefore separating into odd and even entries, where the entries of the index set $I^C_{d,K'}$ containing odd indices are passed through the shifting and scaling network of the coupling blocks. 
The definition of the squeeze and split operations within our network leads to a separation of real and complex components within the split layer. Therefore, our network processes the complex part of the noise signal on a single scale, whereas the real part is processed on two different scales before both signal parts are concatenated in the output layer. Compared to a single-scale network of comparable size we observed a faster convergence behaviour using this architecture, terminating at a number of 25 epochs of training at $b=256$. Moreover, we found an improvement of reconstruction results with respect to both error measures considered in Section~\ref{sec:num} when using the multi-scale architecture in Figure~\ref{fig:network_architecture}.

\begin{figure}[t]
\begin{subfigure}{0.43\textwidth}
\centering
\tikzset{every picture/.style={line width=0.65pt}}   
\begin{tikzpicture}[x=0.65pt,y=0.65pt,yscale=-0.55,xscale=0.55]

\draw  [fill={rgb, 255:red, 255; green, 255; blue, 255 }  ,fill opacity=1 ] (394.49,5) -- (430.65,5) -- (430.65,40.77) -- (394.49,40.77) -- cycle ;

\draw  [fill={rgb, 255:red, 255; green, 255; blue, 255 }  ,fill opacity=1 ] (430.65,5) -- (466.81,5) -- (466.81,40.77) -- (430.65,40.77) -- cycle ;

\draw  [fill={rgb, 255:red, 255; green, 255; blue, 255 }  ,fill opacity=1 ] (467.15,5) -- (503.31,5) -- (503.31,40.77) -- (467.15,40.77) -- cycle ;

\draw  [fill={rgb, 255:red, 255; green, 255; blue, 255 }  ,fill opacity=1 ] (381.74,32.51) -- (417.9,32.51) -- (417.9,68.28) -- (381.74,68.28) -- cycle ;

\draw  [fill={rgb, 255:red, 255; green, 255; blue, 255 }  ,fill opacity=1 ] (417.9,32.51) -- (454.06,32.51) -- (454.06,68.28) -- (417.9,68.28) -- cycle ;

\draw  [fill={rgb, 255:red, 255; green, 255; blue, 255 }  ,fill opacity=1 ] (454.4,32.51) -- (490.56,32.51) -- (490.56,68.28) -- (454.4,68.28) -- cycle ;

\draw  [fill={rgb, 255:red, 255; green, 255; blue, 255 }  ,fill opacity=1 ] (367.84,86.95) -- (404,86.95) -- (404,122.72) -- (367.84,122.72) -- cycle ;

\draw  [fill={rgb, 255:red, 255; green, 255; blue, 255 }  ,fill opacity=1 ] (404,86.95) -- (440.16,86.95) -- (440.16,122.72) -- (404,122.72) -- cycle ;

\draw  [fill={rgb, 255:red, 255; green, 255; blue, 255 }  ,fill opacity=1 ] (440.5,86.95) -- (476.66,86.95) -- (476.66,122.72) -- (440.5,122.72) -- cycle ;

\draw  [fill={rgb, 255:red, 255; green, 255; blue, 255 }  ,fill opacity=1 ] (18.5,80.64) -- (54.66,80.64) -- (54.66,116.41) -- (18.5,116.41) -- cycle ;

\draw  [fill={rgb, 255:red, 255; green, 255; blue, 255 }  ,fill opacity=1 ] (54.66,80.64) -- (90.82,80.64) -- (90.82,116.41) -- (54.66,116.41) -- cycle ;

\draw  [fill={rgb, 255:red, 255; green, 255; blue, 255 }  ,fill opacity=1 ] (90.82,80.64) -- (126.98,80.64) -- (126.98,116.41) -- (90.82,116.41) -- cycle ;
 
\draw  [fill={rgb, 255:red, 255; green, 255; blue, 255 }  ,fill opacity=1 ] (126.98,80.64) -- (163.14,80.64) -- (163.14,116.41) -- (126.98,116.41) -- cycle ;

\draw  [fill={rgb, 255:red, 255; green, 255; blue, 255 }  ,fill opacity=1 ] (199.3,80.64) -- (235.46,80.64) -- (235.46,116.41) -- (199.3,116.41) -- cycle ;

\draw  [fill={rgb, 255:red, 255; green, 255; blue, 255 }  ,fill opacity=1 ] (163.14,80.64) -- (199.3,80.64) -- (199.3,116.41) -- (163.14,116.41) -- cycle ;

\draw  [fill={rgb, 255:red, 255; green, 255; blue, 255 }  ,fill opacity=1 ] (38.78,39.96) -- (74.94,39.96) -- (74.94,75.73) -- (38.78,75.73) -- cycle ;

\draw  [fill={rgb, 255:red, 255; green, 255; blue, 255 }  ,fill opacity=1 ] (74.94,39.96) -- (111.1,39.96) -- (111.1,75.73) -- (74.94,75.73) -- cycle ;

\draw  [fill={rgb, 255:red, 255; green, 255; blue, 255 }  ,fill opacity=1 ] (111.1,39.96) -- (147.26,39.96) -- (147.26,75.73) -- (111.1,75.73) -- cycle ;

\draw  [fill={rgb, 255:red, 255; green, 255; blue, 255 }  ,fill opacity=1 ] (147.26,39.96) -- (183.42,39.96) -- (183.42,75.73) -- (147.26,75.73) -- cycle ;
 
\draw  [fill={rgb, 255:red, 255; green, 255; blue, 255 }  ,fill opacity=1 ] (219.58,39.96) -- (255.74,39.96) -- (255.74,75.73) -- (219.58,75.73) -- cycle ;

\draw  [fill={rgb, 255:red, 255; green, 255; blue, 255 }  ,fill opacity=1 ] (183.42,39.96) -- (219.58,39.96) -- (219.58,75.73) -- (183.42,75.73) -- cycle ;

\draw    (261.99,78.61) -- (350.43,78.72) ;
\draw [shift={(341.43,78.72)}, rotate = 180.09] [fill={rgb, 255:red, 0; green, 0; blue, 0 }  ][line width=0.08]  [draw opacity=0] (-5.93,-4.29) -- (-15,0) -- (-5.93,4.29) -- cycle    ;

\draw  [fill={rgb, 255:red, 255; green, 255; blue, 255 }  ,fill opacity=1 ] (359.15,113.88) -- (395.31,113.88) -- (395.31,149.65) -- (359.15,149.65) -- cycle ;
 
\draw  [fill={rgb, 255:red, 255; green, 255; blue, 255 }  ,fill opacity=1 ] (395.31,113.88) -- (431.47,113.88) -- (431.47,149.65) -- (395.31,149.65) -- cycle ;

\draw  [fill={rgb, 255:red, 255; green, 255; blue, 255 }  ,fill opacity=1 ] (431.81,113.88) -- (467.97,113.88) -- (467.97,149.65) -- (431.81,149.65) -- cycle ;

\draw  [dash pattern={on 4.5pt off 4.5pt}]  (360.65,78.35) -- (550.38,76.32) ;

\draw (260.98,55.11) node [anchor=north west][inner sep=0.75pt]  [font=\small] [align=left] {squeeze};

\draw (500.07,48.24) node [anchor=north west][inner sep=0.75pt]  [font=\normalsize] [align=left] {{\small split}};

\draw (32.12-6,85.48) node [anchor=north west][inner sep=0.75pt]   [align=left] {1};

\draw (67.46-6,85.48) node [anchor=north west][inner sep=0.75pt]   [align=left] {2};

\draw (104.54-6,85.48) node [anchor=north west][inner sep=0.75pt]   [align=left] {3};

\draw (139.87-6,85.05) node [anchor=north west][inner sep=0.75pt]   [align=left] {4};

\draw (176.95-6,85.05) node [anchor=north west][inner sep=0.75pt]   [align=left] {5};

\draw (212.29-6,85.05) node [anchor=north west][inner sep=0.75pt]   [align=left] {6};

\draw (51.82-6,44.79) node [anchor=north west][inner sep=0.75pt]   [align=left] {7};

\draw (87.15-6,44.79) node [anchor=north west][inner sep=0.75pt]   [align=left] {8};

\draw (124.23-6,44.79) node [anchor=north west][inner sep=0.75pt]   [align=left] {9};

\draw (156.65-10,44.37) node [anchor=north west][inner sep=0.75pt]   [align=left] {10};

\draw (193.65-9.5,44.37) node [anchor=north west][inner sep=0.75pt]   [align=left] {11};

\draw (229.07-9.5,44.37) node [anchor=north west][inner sep=0.75pt]   [align=left] {12};

\draw (372.77-6,118.72) node [anchor=north west][inner sep=0.75pt]   [align=left] {1};

\draw (408.11-6,118.72) node [anchor=north west][inner sep=0.75pt]   [align=left] {3};

\draw (444.61-6,118.72) node [anchor=north west][inner sep=0.75pt]   [align=left] {5};

\draw (381.46-6,91.79) node [anchor=north west][inner sep=0.75pt]   [align=left] {2};

\draw (416.8-6,91.79) node [anchor=north west][inner sep=0.75pt]   [align=left] {4};

\draw (453.3-6,91.79) node [anchor=north west][inner sep=0.75pt]   [align=left] {6};

\draw (395.36-6,37.34) node [anchor=north west][inner sep=0.75pt]   [align=left] {7};

\draw (430.7-6,37.34) node [anchor=north west][inner sep=0.75pt]   [align=left] {9};

\draw (464.2-10,37.34) node [anchor=north west][inner sep=0.75pt]   [align=left] {11};

\draw (408.11-6,9.84) node [anchor=north west][inner sep=0.75pt]   [align=left] {8};

\draw (440.53-10,9.84) node [anchor=north west][inner sep=0.75pt]   [align=left] {10};

\draw (477.02-10,9.84) node [anchor=north west][inner sep=0.75pt]   [align=left] {12};
\end{tikzpicture}
\caption{Illustration of a squeeze and a split layer. A squeeze transformation, indicated by the right arrow, transforms an input of shape $p \times K'$ to $2p\times K'/2$. The dashed line indicates a split transformation on the output of the squeeze transformation. This results in two outputs of shape $p\times K'/2$.}
\label{fig:squeeze_split}
\end{subfigure}
\hfill
\begin{subfigure}{0.54\textwidth}
\centering
\tikzset{every picture/.style={line width=0.75pt}}      

\begin{tikzpicture}[x=0.75pt,y=0.75pt,yscale=-0.5,xscale=0.30]

\draw [line width=1.5]    (502,161) -- (946.68,161) ;
\draw [shift={(950.68,161)}, rotate = 179.71] [fill={rgb, 255:red, 0; green, 0; blue, 0 }  ][line width=0.08]  [draw opacity=0] (11.61,-5.58) -- (0,0) -- (11.61,5.58) -- cycle    ;

\draw [line width=1.5]    (503,61) -- (947.68,61) ;
\draw [shift={(951.68,61)}, rotate = 179.71] [fill={rgb, 255:red, 0; green, 0; blue, 0 }  ][line width=0.08]  [draw opacity=0] (11.61,-5.58) -- (0,0) -- (11.61,5.58) -- cycle    ;

\draw [line width=1.5]    (105,88) -- (467.68,88) ;
\draw [shift={(471.68,88)}, rotate = 179.8] [fill={rgb, 255:red, 0; green, 0; blue, 0 }  ][line width=0.08]  [draw opacity=0] (11.61,-5.58) -- (0,0) -- (11.61,5.58) -- cycle    ;

\draw  [color={rgb, 255:red, 0; green, 0; blue, 0 }  ,draw opacity=1 ][fill={rgb, 255:red, 255; green, 255; blue, 255 }  ,fill opacity=1 ][line width=0.75]  (52,2) -- (104.47,2) -- (104.47,175.6) -- (52,175.6) -- cycle ;

\draw [line width=1.5]    (29.87,89.6) -- (47,89.11) ;
\draw [shift={(51,89)}, rotate = 178.37] [fill={rgb, 255:red, 0; green, 0; blue, 0 }  ][line width=0.08]  [draw opacity=0] (11.61,-5.58) -- (0,0) -- (11.61,5.58) -- cycle    ;

\draw  [color={rgb, 255:red, 0; green, 0; blue, 0 }  ,draw opacity=1 ][fill={rgb, 255:red, 255; green, 255; blue, 255 }  ,fill opacity=1 ][line width=0.75]  (133,2) -- (185.47,2) -- (185.47,175.6) -- (133,175.6) -- cycle ;

\draw  [color={rgb, 255:red, 0; green, 0; blue, 0 }  ,draw opacity=1 ][fill={rgb, 255:red, 255; green, 255; blue, 255 }  ,fill opacity=1 ][line width=0.75]  (215,2) -- (267.47,2) -- (267.47,175.6) -- (215,175.6) -- cycle ;

\draw  [color={rgb, 255:red, 0; green, 0; blue, 0 }  ,draw opacity=1 ][fill={rgb, 255:red, 255; green, 255; blue, 255 }  ,fill opacity=1 ][line width=0.75]  (299,1) -- (351.47,1) -- (351.47,174.6) -- (299,174.6) -- cycle ;
 
\draw  [color={rgb, 255:red, 0; green, 0; blue, 0 }  ,draw opacity=1 ][fill={rgb, 255:red, 255; green, 255; blue, 255 }  ,fill opacity=1 ][line width=0.75]  (384,1) -- (436.47,1) -- (436.47,174.6) -- (384,174.6) -- cycle ;

\draw  [color={rgb, 255:red, 0; green, 0; blue, 0 }  ,draw opacity=1 ][fill={rgb, 255:red, 255; green, 255; blue, 255 }  ,fill opacity=1 ][line width=0.75]  (470,1) -- (522.47,1) -- (522.47,174.6) -- (470,174.6) -- cycle ;

\draw  [fill={rgb, 255:red, 255; green, 255; blue, 255 }  ,fill opacity=1 ] (548,1) -- (591.02,1) -- (591.02,143.7) -- (548,143.7) -- cycle ;

\draw  [fill={rgb, 255:red, 255; green, 255; blue, 255 }  ,fill opacity=1 ] (629,1) -- (672.02,1) -- (672.02,143.7) -- (629,143.7) -- cycle ;
 
\draw  [fill={rgb, 255:red, 255; green, 255; blue, 255 }  ,fill opacity=1 ] (715,1) -- (758.02,1) -- (758.02,143.7) -- (715,143.7) -- cycle ;

\draw  [fill={rgb, 255:red, 255; green, 255; blue, 255 }  ,fill opacity=1 ] (797,1) -- (840.02,1) -- (840.02,143.7) -- (797,143.7) -- cycle ;
 
\draw  [fill={rgb, 255:red, 255; green, 255; blue, 255 }  ,fill opacity=1 ] (881,1) -- (924.02,1) -- (924.02,143.7) -- (881,143.7) -- cycle ;

\draw  [color={rgb, 255:red, 0; green, 0; blue, 0 }  ,draw opacity=1 ][fill={rgb, 255:red, 255; green, 255; blue, 255 }  ,fill opacity=1 ][line width=0.75]  (952,1) -- (1004.47,1) -- (1004.47,174.6) -- (952,174.6) -- cycle ;

\draw [line width=1.5]    (1003.87,88) -- (1020.87,88) ;
\draw [shift={(1024.87,88)}, rotate = 182.73] [fill={rgb, 255:red, 0; green, 0; blue, 0 }  ][line width=0.08]  [draw opacity=0] (11.61,-5.58) -- (0,0) -- (11.61,5.58) -- cycle ;

\draw (98.0,9.5) node [anchor=north west][inner sep=0.75pt]  [font=\small,rotate=-90] [align=left] {affine coupling};

\draw (177.5,9.5) node [anchor=north west][inner sep=0.75pt]  [font=\small,rotate=-90] [align=left] {affine coupling};

\draw (93,183.4) node [anchor=north west][inner sep=0.75pt]  [font=\footnotesize, rotate=-90]  {$[ b,p,K']$};

\draw (175,183.4) node [anchor=north west][inner sep=0.75pt]  [font=\footnotesize, rotate=-90]  {$[ b,p,K']$};

\draw (255.5,53.5) node [anchor=north west][inner sep=0.75pt]  [font=\small,rotate=-90] [align=left] {squeeze};

\draw (172+87,183.4) node [anchor=north west][inner sep=0.75pt]  [font=\footnotesize, rotate=-90]  {$[ b,2p,K'/2]$};

\draw (344.5,9.5) node [anchor=north west][inner sep=0.75pt]  [font=\small,rotate=-90] [align=left] {affine coupling};

\draw (-1,79.4) node [anchor=north west][inner sep=0.75pt]  [font=\normalsize]  {$\eta$};

\draw (429.5,9.5) node [anchor=north west][inner sep=0.75pt]  [font=\small,rotate=-90] [align=left] {affine coupling};

\draw (172+2*87,183.4) node [anchor=north west][inner sep=0.75pt]  [font=\footnotesize, rotate=-90]  {$[ b,2p,K'/2]$};

\draw (172+3*87,183.4) node [anchor=north west][inner sep=0.75pt]  [font=\footnotesize, rotate=-90]  {$[ b,2p,K'/2]$};

\draw (515.5,69.5) node [anchor=north west][inner sep=0.75pt]  [font=\small,rotate=-90] [align=left] {split};

\draw (172+4*91,183.4) node [anchor=north west][inner sep=0.75pt]  [font=\footnotesize, rotate=-90]  {$[ b,p,K'/2]$};

\draw (172+4*81,183.4) node [anchor=north west][inner sep=0.75pt]  [font=\footnotesize, rotate=-90]  {$[ b,p,K'/2]$};

\draw (580.5,38.5) node [anchor=north west][inner sep=0.75pt]  [font=\scriptsize,rotate=-90] [align=left] {squeeze};

\draw (587,183.4) node [anchor=north west][inner sep=0.75pt]  [font=\footnotesize, rotate=-90]  {$[ b,2p,K'/4]$};

\draw (665.5,1.5) node [anchor=north west][inner sep=0.75pt]  [font=\scriptsize,rotate=-90] [align=left] {affine coupling};

\draw (582+85,183.4) node [anchor=north west][inner sep=0.75pt]  [font=\footnotesize, rotate=-90]  {$[ b,2p,K'/4]$};

\draw (753.5,1.5) node [anchor=north west][inner sep=0.75pt]  [font=\scriptsize,rotate=-90] [align=left] {affine coupling};

\draw (582+2*86,183.4) node [anchor=north west][inner sep=0.75pt]  [font=\footnotesize, rotate=-90]  {$[ b,2p,K'/4]$};

\draw (833.5,1.5) node [anchor=north west][inner sep=0.75pt]  [font=\scriptsize,rotate=-90] [align=left] {affine coupling};

\draw (582+3*86,183.4) node [anchor=north west][inner sep=0.75pt]  [font=\footnotesize, rotate=-90]  {$[ b,2p,K'/4]$};

\draw (918.5,1.5) node [anchor=north west][inner sep=0.75pt]  [font=\scriptsize,rotate=-90] [align=left] {affine coupling};

\draw (582+4*85,183.4) node [anchor=north west][inner sep=0.75pt]  [font=\footnotesize, rotate=-90]  {$[ b,2p,K'/4]$};

\draw (998.5,10.5) node [anchor=north west][inner sep=0.75pt]  [font=\small,rotate=-90] [align=left] {concatenation};

\draw (1022,80.4) node [anchor=north west][inner sep=0.75pt]  [font=\normalsize]  {$z$};

\draw (995,183.4) node [anchor=north west][inner sep=0.75pt]  [font=\footnotesize, rotate=-90]  {$[ b,p,K']$};
\end{tikzpicture}
\caption{Network architecture of our multi-scale discrepancy network used in Section~\ref{sec:num} with the layer types and corresponding output dimensions for each layer.}
\label{fig:network_architecture}
\end{subfigure}
\end{figure}

\hfill
\begin{table}
    \renewcommand{\arraystretch}{1.1}
    \centering
    \begin{tabular}{|c|c|c|}
    \hline 
     Block & Layer Type & Width \\ \Xhline{3\arrayrulewidth}
     &Fully Connected Layer &  \\ \cline{2-2}
     1 &Layer Normalization & $W_\ell$\\ \cline{2-2}
     & ReLU Activation &  \\ \Xhline{3\arrayrulewidth}
     $\vdots$ & \vdots &$\vdots$ \\ \Xhline{3\arrayrulewidth}
     &Fully Connected Layer &  \\ \cline{2-2}
     6 &Layer Normalization & $W_\ell$\\ \cline{2-2}
     & ReLU Activation &  \\ \Xhline{3\arrayrulewidth}
     \multirow{2}*{7} &Fully Connected Layer & \multirow{2}*{$W_\ell$}\\ \cline{2-2}
     & Tanh Activation & \\ \hline
\end{tabular}
\captionsetup{type=table}

\caption{Network architecture of the scaling and shifting networks $s,t$ of each affine coupling block $\ell = 1,\hdots,8$ utilized in Table~\ref{fig:network_architecture}.}
\label{tab:st_arch}
\end{table}

\vfill

\end{document}